\magnification 1200

\font\tenmsb=msbm10
\font\sevenmsb=msbm10 at 7pt
\font\fivemsb=msbm10 at 5pt
\newfam\msbfam
\textfont\msbfam=\tenmsb
\scriptfont\msbfam=\sevenmsb
\scriptscriptfont\msbfam=\fivemsb
\def\Bbb#1{{\fam\msbfam\relax#1}}

\font\gothicf=eufm10
\font\sgothic=eufm7
\font\ssgothic=eufm5
\textfont5=\gothicf
\scriptfont5=\sgothic
\scriptscriptfont5=\ssgothic
\def\gothic{\fam5}

\font\sc=cmcsc10

\def\gm{{\gothic m}}

\def\ol{\overline}
\def\la{\longrightarrow}
\def\ni{\noindent}
\def\cl{\centerline}

\def\d{\delta }
\def\e{\epsilon}
\def\l{\lambda}
\def\s{\sigma}
\def\a{\alpha}
\def\b{\beta}
\def\g{\gamma}
\def\k{\kappa}
\def\V{V^{d,\d}(\a,\b)}
\def\v{N^{d,\d}(\a,\b)}
\def\D{\Delta}
\def\Dm{\D_m}
\def\Dn{\D_{m-1}}

\def\A{{\Bbb A}}

\def\Q{{\Bbb Q}}
\def\P{{\Bbb P}}

\def\I{{\cal I}}
\def\N{{\cal N}}
\def\O{{\cal O}}
\def\S{{\cal S}}
\def\X{{\cal X}}
\def\Y{{\cal Y}}
\def\U{{\cal U}}

\def\lcm{{\rm lcm}}
\def\th{{\rm th}}
\def\Ker{\mathop{\rm Ker}}
\def\Im{\mathop{\rm Im}}
\def\Coker{\mathop{\rm Coker}}
\def\Spec{\mathop{\rm Spec }}
\def\Div{\mathop{\rm Div }}

\def\Nt{({\cal N } _b)_{{\rm tors}}}

\overfullrule=0pt

\def\sqr#1#2{{\vcenter{\vbox{\hrule height.#2pt \hbox{\vrule width.#2pt height
#1pt \kern #1pt \vrule width.#2pt}\hrule height.#2pt}}}}

\def\qed{{\hfil \break \rightline{$\sqr74$}}}

\

\

\centerline{\bf{COUNTING PLANE CURVES OF ANY GENUS}}

\

\

\

\noindent {\bf Lucia Caporaso}\footnote*{Partially supported by a Sloan
Foundation Fellowship}

\noindent Mathematics department, Harvard University,

\noindent 1 Oxford st., Cambridge MA 02138, USA

\noindent caporaso@abel.harvard.edu

\

\noindent {\bf Joe Harris}

\noindent Mathematics department, Harvard University,

\noindent 1 Oxford st., Cambridge MA 02138, USA

\noindent harris@abel.harvard.edu

\

\

\ni \cl{{\sc Contents}}

\

{\narrower\narrower\smallskip
\noindent 1. Introduction \thinspace \dotfill \thinspace 2

\item{1.} Definitions and notation

\item{2.} The formula for the degrees of generalized Severi varieties

\item{3.} The main results

\item{4.} The formula for irreducible curves

\smallskip

\noindent 2. Geometry of Severi varieties at a general point  \thinspace
\dotfill \thinspace 10

\item{1.} Statement of results

\item{2.} Deformations of maps

\item{3.} Dimension counts and consequences

\item{4.}  Normal sheaves and normal
bundles

\smallskip

\ni 3. Hyperplane sections of Severi varieties: set-theoretic description
\thinspace
\dotfill
\thinspace 23

\item{1.}  The basic setup

\item{2.} Some simplifying assumptions and some corollaries

\item{3.} Proof of simplified Theorem 1.2

\item{4.} The local picture of the degeneration

\item{5.} Verifying the assumptions

\smallskip

\ni 4. Hyperplane sections of Severi varieties: local geometry   \thinspace
\dotfill
\thinspace 33

\item{1.} Deformation spaces of tacnodes

\item{2.} Products of deformation spaces of tacnodes

\item{3.} The local geometry around irreducible curves

\item{4.} The local geometry around reducible curves

\smallskip

\ni References  \thinspace \dotfill \thinspace 47
\smallskip}

\vfill\eject

\

\

\ni \cl{ 1.  {\sc Introduction}}

\

In this paper we  study the geometry of
the {\it Severi varieties} parametrizing plane curves of given degree $d$
and geometric genus
$g$.  As an application, we derive a  recursive formula for their degrees.
This is
the formula in Theorem 1.1,
 which  enumerates the number of curves  in $\P ^2$ of degree $d$ with $\d$
nodes that pass through the appropriate number of points;  that result is
actually
more general, as we shall
explain below.

Such a classical enumerative question has been object of study by many
other people. In
1989,  Z. Ran described a
 different inductive procedure to approach the problem   (cf. [R]). More
recently, interest
in it has been revived by work on quantum cohomology and by the discovery by
M. Kontsevich in 1993 of a beautiful  recursion solving the problem for
curves of genus
$0$ (cf. [KM]; another proof, using different techniques, was given
independently by Ruan and Tian
in [RT]).

Our  approach is  simple. We work over the complex numbers throughout. We
denote by $\P^N$ the projective space
of all plane curves of degree $d$ and by $V^{d,\d} \subset \P^N$ the
closure of the
subset of $\P^N$ corresponding to curves having exactly $\d$ nodes as
singularities.
Also, for any point $p \in \P^2$ we  let $H_p \subset\P^N$ be the hyperplane of
curves containing the point $p$. Our procedure consists in intersecting the
variety
$V^{d,\d}$ with a succession of hyperplanes of the form
$H_{p_i}$, where the points $p_i$ are general points on a fixed line $L
\subset \P^2$. At
each stage we are able to describe the irreducible components of the
intersection;
the point is, they all belong to a specific collection of varieties,
which we call {\it generalized Severi varieties} and denote by
$\V$. These pararmetrize plane curves of given
degree $d' \le d$ and genus $g' \le g$ satisfying certain tangency
conditions with respect
to the line $L$. More generally, we can express the intersection of any
generalized Severi
variety $\V$ with a hyperplane $H_p$ corresponding to a general point $p
\in L$ as a
union of generalized Severi varieties
$V^{d',\d'}(\a',\b')$ of dimension one less; counting multiplicities
correctly, this allows us
to derive our recursive statement.

	\

\noindent  1.1. {\bf  Notation and definitions.}
We now introduce the notation and precise definitions that will allow us
to state our formula.

For any sequence $\a =
(\a_1,\a_2,\ldots)$ of nonnegative integers with all but finitely many
$\a_i$ zero, set
$$
\#\a \; = \; \#\{i : \a_i \ne 0\}
$$
$$ |\a| \; = \; \a_1+\a_2+\dots+\a_n
$$
$$ I\a \; = \; \a_1+2\a_2+\dots+n\a_n
$$ and
$$ I^\a \; = \; 1^{\a_1}2^{\a_2}3^{\a_3}\cdots .
$$
We   denote by $\lcm(\a)$ the least common multiple of the set $\#\{i : \a_i \ne
0\}$.

We  denote by
$e_k$ the sequence
$(0,\ldots,0,1,0,\ldots)$ that is zero except for a $1$ in the $k^{\rm th}$
term (so that any
sequence $\a =
 (\a_1, \a_2, \ldots)$ is expressible as $\a = \sum \a_ke_k$). By the inequality
$\a
\ge
\a'$ we  mean
$\a_k
\ge
\a'_k$ for all
$k$; for such a pair of sequences we  set
$$ {\a \choose \a'} \; = \; {\a_1 \choose \a'_1}{\a_2 \choose \a'_2}{\a_3
\choose
\a'_3}
\cdots
$$

We now define the main objects of study, the
varieties $\V$ parametrizing plane curves of given degree and geometric genus
satisfying certain tangency conditions with respect to a line. Fix a line
$L
\subset
\P^2$ and a collection
$$
\Omega \; = \; \{p_{i,j}\}_{1\le j \le \a_i} \subset L
$$
 of $|\a|$  general
points on $L$. For any
$d, \d,
\a$ and
$\b$ satisfying
$I\a+I\b=d$, we define the {\it generalized Severi variety}
$\V(\Omega)$ to be the closure of the locus of reduced plane curves $X$ of
degree
$d$ and geometric genus $g = {d-1 \choose 2} - \d$, not containing $L$,  with
(informally)
$\a_k$ ``assigned" points of contact of order $k$ and $\b_k$ ``unassigned"
points of
contact of order $k$ with $L$. Formally,  we require that, if
$\nu : X^\nu \to X$ is the normalization of $X$, then there exist $|\a|$ points
$q_{i,j}
\in X^\nu$, $j = 1,\ldots,\a_i$, and $|\b|$ points
$r_{i,j}
\in X^\nu$, $j = 1,\ldots,\b_i$, such that
$$
\nu(q_{i,j}) \; = \; p_{i,j}
$$ and
$$
\nu^*(L) \; = \; \sum i\cdot q_{i,j} \; + \; \sum  i\cdot r_{i,j} .
$$
Where the dependency on the points $p_{i,j}$ is not relevant---for example, in
discussions of the dimension or degrees of generalized Severi
varieties---we will often
suppress the $\Omega$.

To take some simple cases, taking $\a=0$ and $\b=(d,0,\ldots)$ imposes no
condition at all, that is, $V^{d,\d}((0,0,\ldots),(d,0,\ldots))$ is simply
the closure
$V^{d,\d}$ of the locus of plane curves of degree
$d$ with $\d$ nodes. Taking $\a=(1,0,\ldots)$ and $\b=(d-1,0,\ldots)$ we
get the closure
of the locus of such curves passing through a single fixed point of $L$;
and taking
$\a=0$ and
$\b=(d-2,1,0,\ldots)$ we get the closure of the locus of such curves
tangent to $L$ at a
smooth point of the curve.

Note that we do not require
$X$ to be irreducible. Classically, the term ``Severi variety" means a variety
parametrizing irreducible curves of given degree and genus, so we are
somewhat at odds
with traditional usage here; but we will find it much more convenient, in
both the
statement and proof of the results below, to include components of $\V$
whose general
member is reducible.

Let $V$ be a possibly reducible variety. We will say that $V$ has pure
dimension if
all irreducible components of $V$ have the same dimension. Moreover
whenever  we make a statement about the general point
of $V$, we mean that the statement holds for a general point of
any irreducible component of $V$.

We will adopt the following convention, we will   denote the various types
of Severi varieties by the
 symbol ``$V$" to which we will add certain decorations, and we will
correspondingly
use the symbol ``$N$" with same
decorations to denote the degree of $V$ as a subvariety of $\P ^N$; for
example, we define
$\v := \deg \V$.

\

\noindent 1.2. {\bf The formula for the degrees of generalized Severi
varieties.}
As a
result of the analysis of hyperplane sections of generalized Severi varieties
 we have the recursive formula

\proclaim Theorem 1.1.  Let
$\v$ be the degree of $\V$. Then
$$
\eqalign{\v \; &= \; \sum_{k : \b_k > 0} k \cdot N^{d,\d}(\a+e_k,\b-e_k)
\cr &\quad + \sum I^{\b'-\b}{\a \choose \a'}{\b' \choose
\b}N^{d-1,\d'}(\a',\b') \cr}
$$ where the second sum is taken over all $\a', \b'$ and $\d' \ge 0$ satisfying
$$
\eqalign{\a' &\le \a \cr
\b' &\ge \b \cr
\d' &\le \d \cr
\d-\d' + |\b'-\b| &= d-1  \cr}
$$

Taking $\a=0$ and $\b = (d,0,\ldots)$, we get the degree of the closure of
the variety of (not necessarily irreducible) plane curves of degree
$d$ with $\d$ nodes. We can find the degree of the component parametrizing
irreducible such curves---that is, the classical Severi variety---by
subtracting off the
degrees of the others, which we know recursively. Alternatively, we can
give a recursion
formula to calculate directly the degrees of the varieties parametrizing
irreducible plane
curves, and will do so in the last section of this chapter; but this
formula is more
complicated.

We  now illustrate how the above formula works by computing the degree of
the Severi variety
of quartics with three nodes (we  assume known the degrees of the
generalized Severi
varieties parametrizing cubics satisfying tangency conditions).  To shorten
notation, we
write
$(d,\d,\a,\b)$ for
$N = N^{d,\d}(\a,\b)$, and suppress the zeroes at the end of sequences $\a$
and $\b$ and the
parentheses around sequences
$\a$ and
$\b$ of length 1. The result of intersecting the variety $V^{4,3}(0,4)$
with five successive
hyperplanes of the form $H_p$ is then the following five equations  (we
denote the
contribution of each component to the degree, where known, in angle brackets).
$$
\eqalign{(4,3,0,4) \; = \; (4,3,1,3) \; &= \; (4,3,2,2) \cr &\quad + \;
(3,0,0,3) \quad
\langle \underline{1} \rangle\cr}
$$
$$
\eqalign{(4,3,2,2) \; &= \; (4,3,3,1) \cr  &\quad + \; 3(3,1,0,3) \quad
\langle 3 \times
12 =
\underline{36} \rangle\cr  &\quad + \; 2(3,0,1,2) \quad \langle 2 \times 1 =
\underline{2} \rangle\cr}
$$
$$
\eqalign{(4,3,3,1) \; &= \; (4,3,4,0) \cr  &\quad + \; 3(3,2,0,3) \quad
\langle 3 \times
21 =
\underline{63} \rangle\cr  &\quad + \; 2(3,1,0,(1,1)) \quad \langle 2
\times 36 =
\underline{72} \rangle\cr &\quad + \; 6(3,1,1,2) \quad \langle 6 \times 12 =
\underline{72} \rangle\cr &\quad + \; 3(3,0,2,1) \quad \langle 3 \times 1 =
\underline{3} \rangle\cr}
$$ and finally
$$
\eqalign{(4,3,4,0) \; &= \; (3,3,0,3) \quad \langle
\underline{15} \rangle\cr  &\quad + \; 4(3,2,1,3) \quad \langle 4 \times 21 =
\underline{84} \rangle\cr &\quad + \; 2(3,2,0,(1,1)) \quad \langle 2 \times 30 =
\underline{60} \rangle\cr &\quad + \; 6(3,1,2,1) \quad \langle 6 \times 12 =
\underline{72} \rangle \cr &\quad + \; 8(3,1,1,(0,1)) \quad \langle 8
\times 16 =
\underline{128} \rangle \cr &\quad + \; 3(3,1,0,(0,0,1)) \quad \langle 3 \times 21 =
\underline{63} \rangle \cr &\quad + \; 4(3,0,3,0) \quad \langle 4 \times 1 =
\underline{4} \rangle\cr} .
$$ Adding it all up, we find that
$$
\eqalign{(4,3,4,0) \; &= \; 15 + 84+ 60 + 72 + 128 + 63 + 4  \; = \;
426\cr (4,3,3,1) \; &=
\; 426 + 63 + 72 + 72 + 3  \; = \;  636\cr (4,3,2,2) \; &= \; 636 + 36 + 2
\; = \;  674\cr}
$$ and so
$$ (4,3,0,4) \; = \; (4,3,1,3) \; = \; 674 + 1 \; = \; 675 .
$$

Now, $V^{4,3}((0),(4))$  has two irreducible components  of dimension 11,
one coincides withe the classical $V_{4,3}$ which  parametrizes irreducible
quartics
with three nodes; the other parametrizes reducible curves that are the union
of a line and a cubic and has degrre ${11 \choose 2} = 55$.
and so we conclude that the degree of the classical Severi variety is
$620$.

R. Vakil has checked the formula in all degrees up to and
including 6; and the results agree with those of Vainsencher for $\d \le 6$
and with
those of Kontsevich-Manin for
$g = 0$.  We  list below the numbers obtained by R. Vakil (by applying this
formula)
for the degrees $N$ of the Severi varieties $V^{d,\d}(0,d)$
 for $d = 5$ and 6 and all
possible values of $\d$.

\

\settabs \+ \indent  \indent & \qquad \qquad&\qquad\qquad & \qquad\qquad  &
59809860
& \cr

\+&$d$&$\d$&$g$&\hfill $N$ & \cr
\bigskip
\+&5&0&6&\hfill1& \cr
\+&&1&5&\hfill48& \cr
\+&&2&4&\hfill882& \cr
\+&&3&3&\hfill7915& \cr
\+&&4&2&\hfill36975& \cr
\+&&5&1&\hfill90027& \cr
\+&&6&0&\hfill109781& \cr
\+&&7&-1&\hfill65949& \cr
\+&&8&-2&\hfill26136& \cr
\+&&9&-3&\hfill6930& \cr
\+&&10&-4&\hfill945& \cr
\medskip
\+&6&0&10&\hfill1& \cr
\+&&1&9&\hfill75& \cr
\+&&2&8&\hfill2370& \cr
\+&&3&7&\hfill41310& \cr
\+&&4&6&\hfill437517& \cr
\+&&5&5&\hfill2931831& \cr
\+&&6&4&\hfill12597900& \cr
\+&&7&3&\hfill34602705& \cr
\+&&8&2&\hfill59809860& \cr
\+&&9&1&\hfill63338881& \cr
\+&&10&0&\hfill40047888& \cr
\+&&11&-1&\hfill15580020& \cr
\+&&12&-2&\hfill4361721& \cr
\+&&13&-3&\hfill718918& \cr
\+&&14&-4&\hfill135135& \cr
\+&&15&-5&\hfill10395& \cr

\

\noindent  1.3. {\bf The main results.}
As we indicated at the outset, what we
actually prove is an equality of cycles rather than of numbers: we  show
first that
for $p \in L$ general the intersection
$\V
\cap H_p$ is a union of varieties of the form
$V^{d',\d'}(\a',\b')$, and say which ones occur; and then we will calculate
the intersection multiplicity of $\V$
and $H_p$ along each such variety.

\

\ni \underbar{Remark}.
Throughout, we will identify the projective space of plane curves of degree
$d-1$
with the subspace of plane curves of degree $d$ containing $L$. Thus, for
example, by
$V^{d-1,\d'}(\a',\b')(\Omega')$ we mean the closure in $\P^N$ of the locus
of curves
$X_0 = X \cup L$ where $X$ is a plane curve of degree $d-1$, not containing
$L$, having
$\d'$ nodes and $\a'$ assigned and $\b'$ unassigned points of contact with $L$.

\proclaim Theorem 1.2.  The intersection $\V(\Omega)
\cap H_p$ is contained in a union of varieties (without common components)
as follows:
\item {a.} For each $k$ such that $\b_k
> 0$, the variety
$$
V^{d,\d}(\a+e_k,\b-e_k)(\Omega \cup \{p_{k,\a_k+1}=p\}) \, ;
$$
and
\item {b.} For each $\a' \le \a$, $\b' \ge \b$ and $\d' \le \d$ with
$\d-\d' + |\b'-\b|
= d-1$, the union of the varieties
$$V^{d-1,\d'}(\a',\b')(\Omega')$$
where $\Omega'$ ranges over all subsets $\Omega' = \{p_{i,j}'\}_{1\le j \le
\a'_i} \subset
\Omega$ such that $\{p'_{i,1},\ldots,p'_{i,\a'_i}\} \subset
\{p_{i,1},\ldots,p_{i,\a_i}\}$ for
each $i$.

\ni \underbar{Remarks}.
1.
We can express the last condition $\d-\d' + |\b'-\b|
= d-1$ on $\d'$ and $\b'$  in terms of the geometric genera $g = {d-1
\choose 2} - \d$ and $g' = {d-2
\choose 2} - \d'$: we have
$$
\eqalign{g - g' \; &= \; ({d-1
\choose 2} - \d) - ({d-2
\choose 2} - \d') \cr
&= \; d-2 - (\d - \d') \cr
&= \; |\b' - \b| - 1 \, . \cr}
$$

2. Note that there are a total of ${\a \choose \a'}$ varieties of
the form $V^{d-1,\d'}(\a',\b')$ in part (b) of this statement, which
accounts for the factor
${\a \choose \a'}$ in the formula in Theorem 1.1. The remaining factors will be
intersection multiplicities, as described in Theorem 1.3 below.

3. By the  dimension counts of Section 2, all the
varieties listed in the statement of Theorem 1.2 have pure dimension
$\dim(\V)-1$; so it follows that the intersection $\V(\Omega) \cap H_p$
will consist of
the union of a subset of these. In fact the intersection is equal to the
union of all of them, as will follow from the analysis of the local
geometry of $\V$
given in the proof of Theorem 1.3.

Having described the
intersection set-theoretically, we  now ask about the local geometry of the
larger
variety $\V$ along each component of the intersection: how many branches it
has, and
with what multiplicity each intersects the hyperplane
$H_p$. The answer to both (and hence the multiplicity with which each
component of $\V
\cap H_p$ appears in the intersection cycle) is expressed in the following.

\proclaim Theorem 1.3.
\item {a.}  Let $V'=V^{d,\d}(\a+e_k,\b-e_k)(\Omega \cup \{p\})$ be as in
part (a)
of Theorem 1.2. Then $V' \subset \V(\Omega) \cap H_p$, and at a general point of
$V'$ the variety $\V$ is smooth and has
intersection multiplicity $k$ with $H_p$ along $V'$.
\item  {b.} Let $V' = V^{d-1,\d'}(\a',\b')(\Omega')$ be as in part (b) of
Theorem 1.2. At a
general point of
$V'$, the variety $\V$ will have ${\b' \choose \b}I^{\b'-\b}/\lcm(\b' -
\b)$  branches,
each of which will  have intersection multiplicity $\lcm(\b'-\b)$ with
$H_p$ along $V'$.

The recursive formula for the degrees $\v$ of the generalized Severi
varieties given in
Theorem 1.1 follows directly from these two statements.

The proofs of Theorems 1.2 and 1.3 will be given in Chapters 3 and 4,
following some
preliminary deformation-theoretic arguments and dimension counts in Section
2. For the
most part the proof of Theorem 1.2 follows the lines of that of Proposition
2.5  of
[CH]. To prove it, we will use the technique of semistable reduction to
analyze a
family of curves
$X
\in \V$ specializing to a curve
$X_0$ containing $L$.  This approach yields a number conditions that the
curve $X_0$
corresponding to a general point of the intersection $\V
\cap H_p$ must satisfy, some of which are far from obvious
from the point of view of the geometry of plane curves alone. We then
compare these
with the dimension estimates of Chapter 2, using the fact that $X_0$ is a
general
member of a family of dimension
$\dim(\V)-1$, to obtain  an exact description of the set-theoretic
intersection $\V \cap H_p$.

As for Theorem 1.3, this requires a deeper analysis of the local structure
of the total space of such a family, based on the deformation theory of the
tacnodes of
$X_0$.  As in the case of Theorem 1.2, this is based on arguments in [CH];
but while the
proof of Theorem 1.2 is largely parallel to the corresponding argument of
[CH], the
argument for Theorem 1.3 requires additional work. Briefly, the proof of
Proposition
2.7 in [CH] rests on the description given there of the deformation space
of a single
tacnode; this suffices for the purposes of that paper.
Here we do need to
consider degenerations having more than one tacnode, hence we develop
an analysis of the geometry of a product of deformation spaces of tacnodes,
building on the
description given in [CH] of the deformation space of a single tacnode.
We should mention that some of the results
on deformations of tacnodes have also been obtained by Z. Ran in [R].

\

\ni  1.4. {\bf The formula for irreducible curves.}
We now  give a formula for the degrees of the varieties parametrizing
irreducible plane curves of given degree and genus satisfying tangency
conditions.

Denote by $V_{d,\d}(\a,\b)$ the union of the components of $\V$
whose general point $[X]$ corresponds to an irreducible curve $X \subset
\P^2$, and by
$N_{d,\d}(\a,\b)$ its degree. Now, consider the intersection of a variety
$V_{d,\d}(\a,\b)$ with
the hyperplane $H_p$, with $p \in L$, and let $[X_0]$ be a general
point  of a component $V$ of the intersection. If $X_0 = X \cup L$ and $X$
has irreducible
components
$X_1,\ldots,X_k$ of degrees $d_1,\ldots,d_k$, then each component $X_j$
will correspond to a
general point of a variety $V_{d_j,\d_j}(\a^j,\b^j)$. The part of the
intersection
$V_{d,\d}(\a,\b) \cap H_p$ corresponding to curves containing $L$ will thus
be a Segre image
of a product of varieties $V_{d_j,\d_j}(\a^j,\b^j)$, and its degree will be
the product of the
degrees of the factors $V_{d_j,\d_j}(\a^j,\b^j)$, times a multinomial
coming from the formula
for the degrees of Segre images. We arrive in this way at the formula
$$
\eqalign{N_{d,\d}(\a,\b) \; &= \; \sum_{k : \b_k > 0} k \cdot
N_{d,\d}(\a+e_k,\b-e_k)
\cr &\quad + \sum {1 \over \s}
{2d+g-2+|\b| \choose 2d_1+g_1-1+|\b^1|,\ldots, 2d_k+g_k-1+|\b^k|} \cdot
{\a \choose \a^1,\ldots,\a^k} \cdot \cr
&\qquad \qquad \cdot \prod_{j=1}^k{\b^j+\g^j \choose \b^j} \cdot
\prod_{j=1}^kI^{\g^j}
\cdot \prod_{j=1}^k N_{d_j,\d_j}(\a^j,\b^j + \g^j)\cr}
$$
where the second sum is taken over all collections of integers
$d_1,\ldots,d_k$ and
$\d_1,\ldots,\d_k$ and collections of sequences
$\a^1,\ldots,\a^k$, $\b^1,\ldots,\b^k$ and $\g^1,\ldots,\g^k$
satisfying
$$
\eqalign{\a^1+\ldots +\a^k \; &\le  \; \a \cr
\b^1+\ldots +\b^k  \; &= \;  \b \cr
|\g^j|  \; &>  \; 0 \cr
d_1+\ldots +d_k \;  &= \; d-1  \quad {\rm and} \cr
\d_1+\ldots+\d_k \; &= \;  \d + \sum|\g^j| - \sum_{i<j} d_id_j - d + 1 \, .\cr}
$$
Here  by the symbol ${n \choose a_1,\ldots,a_k}$ we mean the multinomial
coefficient
$$
{n \choose a_1,\ldots,a_k} \; = \; {n\! \over a_1!\cdots a_k! (n-a_1-\ldots
-a_k)!}
$$
and correspondingly for a collection of sequences $\a$ and
$\a^1,\ldots,\a^k$ we set
$$
{\a \choose \a^1,\ldots,\a^k} \; = \; \prod_i {\a_i \choose
\a_i^1,\ldots,\a_i^k} \, .
$$
By $g_j$ we mean ${d_j-1 \choose 2} - \d_j$.

The symbol $\s$ is 1 except in rare cases. It is the degree of the map from
the union of
the product of varieties of the form
$V_{d_j,\d_j}(\a^j,\b^j)$ to its image in $\P^N$: given the integers
$d_1,\ldots,d_k$ and
$\d_1,\ldots,\d_k$ and collections of sequences
$\a^1,\ldots,\a^k$, $\b^1,\ldots,\b^k$ and $\g^1,\ldots,\g^k$ we define an
equivalence relation
on the set $\{1,2,\ldots,k\}$ by saying $i \sim j$ if $d_i = d_j$,
$\d_i = \d_j$, $\a^i=\a^j$, $\b^i=\b^j$ and $\g^i=\g^j$ and define $\s$ to
be the product of the
factorials of the cardinalities of the equivalence classes.

This formula follows from Theorems 1.2 and 1.3  in much the same way as
Theorem 1.1. Note
that we have here decomposed $\b$ into $\sum \b^j$ and the difference $\b'
- \b$ into $\sum
\g^j$, and further specified that $|\g^j| > 0$. This is because (as we will
see in Chapter 3) the
``new" unassigned points  of $X\cap L$
(that is, the points of $X \cap L$ that are not limits of points of
intersection of nearby
curves in $V_{d,\d}(\a,\b)$ with $L$) correspond to points of intersection
of the
normalizations of the components $X_j$ with $L$ in the nodal reduction of
the family. Since
we are only concerned here with curves $X$ arising as limits in families of
irreducible
curves, their nodal reductions must be connected; and this corresponds to
the requirement
$|\g^j| > 0$.

\vfill\eject

\ni \cl{{\sc 2. Geometry of Severi varieties at a general point}}

\

\ni  2.1. {\bf Statement of results.}
In this section we will compute the dimensions of generalized Severi
varieties $V$ and we will describe the geometry of its general points.

A naive reasoning yields a lower bound for the dimension of $V$. Namely,
if we impose
no conditions on the intersections of our curves with the line $L$, the
corresponding
locus---the classical Severi variety--- has codimension $\d$ in the space
$\P^N$ .
 Requiring that a curve $X$ have
intersection multiplicity $i$ with $L$ at a specified point $p_{i,j}$ is
$i$ linear conditions on the coefficients of $X$,
which we would  expect to be independent; and if we don't specify the point, the
codimension of the corresponding locus should be one less, that is, $i-1$.
In sum we have
$$
\eqalign{\dim(\V) \; &\geq \; {d+2 \choose 2} - 1 - \d - I\a - (I\b-|\b|) \cr
&= \;  {d+1 \choose 2} - \d + |\b| \, . \cr}
$$
Or,
   in terms of the geometric
genus of the curves involved,
$$
\eqalign{\dim(\V) \; &\geq  \; 3d + g - 1  - I\a - (I\b-|\b|) \cr
&= \;  2d + g - 1 + |\b| \, . \cr}
$$

We shall prove that equality holds:

\proclaim Proposition 2.1.   $\V$ has pure dimension $2d +
g - 1 + |\b|$.

Likewise, there are no surprises when it comes to the geometry of a general
member
$X$ of a generalized Severi variety $\V$. We would expect
the curve $X$ to have only nodes as singularities, to be smooth at its points of
intersection with $L$, and so on; and this is indeed the case. We list the
relevant facts  in the following Proposition.

To do that,  fix any curve
$G \subset \P^2$ and any finite subset $\Gamma \subset \P^2$. Let $[X] \in
\V$ be a
general point of a generalized Severi variety, $X \subset \P^2$ the
corresponding curve,
$\nu : X^\nu \to X \subset \P^2$ its normalization. Let $\{q_{i,j}\}$ and
$\{r_{i,j}\} \subset
X^\nu$ be such that $\nu(q_{i,j}) = p_{i,j}$ and
$$
\nu^*L \; = \; \sum i \cdot q_{i,j} + \sum i \cdot r_{i,j}
$$
and let $s_{i,j} = \nu(r_{i,j}) \in L$ be the image of $r_{i,j}$ (that is,
$\{s_{i,j}\} \subset
L$ will be the ``unassigned points" of intersection of $X$ with $L$). We
have then the

\proclaim Proposition 2.2.
\item{a.} $X$ has only nodes as singularities.
\item{b.}  $X$ is smooth along $X \cap L$.
\item{c.}  The points $\{q_{i,j}\}$ and $\{r_{i,j}\} \subset
X^\nu$ are all distinct.
\item{d.}  The points $\{p_{i,j}\}$ and $\{s_{i,j}\} \subset
L$ are all distinct.
\item{e.} $X$ intersects $G$ transversely (in particular, $X$ is smooth
along $X \cap G$)
and is disjoint from $\Gamma$.

To prove these statements, we will need some results  about deformations of
maps, which will be the object of the next section.

\

\ni  2.2. {\bf Deformations of maps.}
Throughout, we will assume we are working over a field of
characteristic zero, and will use the analytic topology where necessary.

We will be concerned with families of maps from a possibly variable smooth
domain to a fixed
smooth target space. In other words, we will consider a flat, smooth,
proper family $f : \X \to B$
over a  smooth connected base $B$, a smooth variety $\Y$ and a morphism
$\psi :
\X
\to B
\times Y$ of
$B$-schemes. For each
$b \in B$, we let $\psi_b : X_b \to Y$ be the restriction of $\psi$ to the
fiber $X_b$ of $\X$ over
$b$, and
$$
d\psi_b \; : \; TX_b \; \la \; \psi_b^*TY
$$
the differential of $\psi_b$.
We let $\N_b$ be the normal sheaf of $\psi_b$, that is, the cokernel of the
morphism $d\psi_b$
of sheaves on $X_b$. Equivalently, if we let
$$
d\psi \; : \; T\X \; \la \; \psi^*T(B \times Y)
$$
be the differential of $\psi$ and $\N = \Coker(d\psi)$ the normal sheaf of
$\psi$, then the
normal sheaf $\N_b$ of $\psi_b$  is the restriction of $\N$ to the fiber
$X_b$, that is,
$
\N_b \; = \; \N \otimes \O_{X_b} \, .
$
Note that if $\psi_b$ is an immersion, then $\N_b$ will be locally free;
more generally, if
$\psi_b$ is equidimensional onto its image then the sheaf
$\N_b$ will have a torsion subsheaf supported exactly on the locus where
$d\psi_b$ fails to be
an injective bundle map.

We  now describe  the {\it Kodaira-Spencer map}  of the family $\psi$ of
morphisms. This is a map
$\k : T_bB \la H^0(X_b, \N_b)$
that associates to any tangent vector $v \in T_bB$ to $B$ a global
section $\s = \k(v)$ of the normal sheaf,  in such a way that the family is
trivial (that is, the
family $\X \cong B \times X_b$ as $B$-schemes and the morphism $\psi = id_B
\times
\psi_b$---if and only if $\k(v) = 0$ for every $v$).
To define it, let
 $\pi : B \times Y \to B$ be the projection, we have an inclusion of
bundles
$$
\pi^*T_B \; \hookrightarrow \; T_{B \times Y} \, .
$$
 We let $i :
\psi^*\pi ^*TB \hookrightarrow \psi^*T(B \times Y)$ be the corresponding
inclusion of pullbacks to
$\X$, and let
$\tilde \k :  \psi^*\pi ^*TB\to \N$ be the composition of $i$ with the
 surjection $\psi^*T_{B \times Y} \to \N$.

\

\vskip1.4in

\hskip.5in \special {picture KSdiagram}

\

\ni Restricting to $X_b$ and taking global sections, we get a map
$$
\k_b \; : \; T_bB \; \hookrightarrow \; H^0(X_b, \psi^*\pi^*T_B)  \; \la \;
H^0(X_b,\N_b)
$$
which we will call the {\it Kodaira-Spencer map} of the given family at $b$.
Equivalently, we let $\k$ be the pushforward of $\tilde \k$ to $B$,
composed with the inclusion
of $T_B$ into $ f_*\psi^*\pi^*T_B$: that is,
$$
\k \; = \; f_*\tilde \k \; : \; T_B \; \hookrightarrow \; f_*\psi^*\pi^*T_B
\; \la \; f_*\N \, .
$$
We will call $\k$ the global Kodaira-Spencer map of the family; the maps
$\k_b$ are then
 the composition of the induced maps $T_bB \to (f_*\N)_b$  on stalks with
the natural
maps
$(f_*\N)_b \to H^0(X_b,N_b)$.

The standard applications of this construction rest on two facts. The first
is that if the family
$\psi$ of morphisms is nowhere isotrivial (that is, the restriction  of
$\psi$ to the subfamily
$\X_{B_0} = f^{-1}(B_0) \subset \X$ is not trivial for any analytic arc $B_0
\subset B$), then at a general point
$b
\in B$ the map
$\k_b$ must be injective, so that we have an a priori bound on the
dimension of the family:
$$
\dim(B) \; \le \; h^0(X_b, \N_b) \, .
$$
(If, for general $b \in B$, we had
$\Ker(\kappa_b) \ne 0$, we could in an
analytic neighborhood of $b$ restrict to a curve whose tangent space was
contained in
$\Ker(\kappa_b)$ at each point.) Secondly, the chern classes of the normal
sheaf are in general
readily calculated, so that in many cases it may be possible to estimate
$h^0(X_b, \N_b)$.

In the
case of plane curves, for example, if $X_b$ is a curve of genus $g$ and
$\psi_b : X_b \to \P^2$ is
birational onto a plane curve of degree $d$, then $\N_b$ is a rank one
sheaf on the curve
$X_b$, the degree of whose chern class is
$$
\eqalign{\deg(c_1(\N_b)) \; &= \; \deg(c_1(\psi_b^*T_{\P^2})) -
\deg(c_1(T_{X_b})) \cr
&= \; 3d + 2g - 2 \cr
&> \; 2g - 2 \, .\cr }
$$
We would thus
expect that
$$
\eqalign{\dim(B)  \; &\le \; h^0(X_b, \N_b) \cr
&= \;  \deg(c_1(\N_b)) - g + 1 \cr
&= \; 3d + g - 1 \, . \cr}
$$
We cannot, however, conclude this yet. The difficulty arises from the
possibility that $\psi_b$ is
not an immersion: if the differential
$d\psi_b$ vanishes at points of $X_b$, the sheaf $\N_b$ will have torsion
there, and in this case
the quotient $\N_b/(\N_b)_{{\rm tors}}$ (and hence $\N_b$ itself) may well
be special. In such a
case, the dimension $h^0(X_b, \N_b)$ will indeed be larger than the naive
estimate $3d+g-1$ for
the dimension of our family, and the method appears to fail.

Happily, there is a standard result  that deals with this situation.
The current version was worked out in conversations
with Johan de Jong, to whom we are very grateful.

Let  $\X \to B$  be as before and
assume that  $\psi : \X \to B \times Y$  is
birational onto its image.
Then we have

\proclaim Lemma 2.3.  If $b \in B$ is a general point, then
$$
\Im(\kappa_b) \; \cap \; H^0(X_b, (\N_b)_{{\rm tors}}) \; = \; 0 \, .
$$

\ni \underbar{Remarks}.
1. If we do not assume the map $\psi$ is birational onto its image, the
conclusion of the
Lemma may well be false. In fact, it will fail exactly when the map $\psi_b
: X_b \to Y$ is
multiple-to-one, with constant image but variable branch points.

2. While we will not introduce the definitions needed to make
this precise,  another way to express this Lemma is to say that ``the
first-order deformation of the map
$\psi_b$ corresponding to a torsion section of $N_b$ can never be
equisingular". If $b \in B$ is
general the first-order deformations of $\psi_b$ arising from the family
$\psi : \X \to B \times Y$
are necessarily equisingular; it follows that they cannot be torsion.

\

\ni {\it Proof}. Note first that, using the analytic topology, it is enough to
prove the Lemma in case $B$ is one-dimensional:  if we had
$\Im(\kappa_b) \cap  H^0(X_b, (\Nt ) \ne 0$ at general $b \in B$ we could in an
analytic neighborhood of $b$ restrict to a curve whose tangent space was
contained in
$(\kappa_b)^{-1}(H^0(X_b, (\N_b)_{{\rm tors}}))$ at each point.

We may thus assume that $\psi : \X \to B \times Y$ is a one-parameter
family of maps, the
image of whose Kodaira-Spencer map $\k_b$ at a general point is contained in $H^0(X_b, \Nt )$.
Let $Z = \psi(\X) \subset B \times Y$ be the image of $\X$,  $p \in \X$ a
general
point with image $\psi(p) = (b,q) \in B \times Y$. We are assuming that for any
 $v \in T_bB$, the image $\k_b(v)$ vanishes at $p$; that is, the tangent space
$T_{(b,q)}Z$ is of the form
$$
T_{(b,q)}Z \; = \; T_bB \times \Lambda_p
$$
for some linear subspace $\Lambda_p \subset T_qY$.

Now, let $t$ be a local analytic coordinate on $B$ near $b$, and
$(x,y_1,\ldots,y_n)$ local
coordinates on $Y$ near $q$ such that $\psi_b^*x$ is a local coordinate on
$X_b$ near $p$ (so
that  the pair $(t,x)$ give local coordinates on the surface $\X$ near
$p$). We can write
the map $\psi$ locally as
$$
y_i \; = \; f_i(t,x) , \qquad i = 1,\ldots,n \, .
$$
The tangent space $T_{(b,q)}Z$ is then the zero locus of the linear forms
$$
dy_i - {\partial f_i \over \partial t}dt - {\partial f_i \over \partial x} dx
$$
and that statement that $T_{(b,q)}Z \; = \; T_bB \times \Lambda_p
$ for some linear subspace $\Lambda_p \subset T_qY$ says that ${\partial
f_i \over \partial t}$
vanishes identically near $p$. We deduce that the image of
$\psi_b$ is constant, i.e., that near $(b,q)$ the image $Z$ is equal to the
product of a
neighborhood of $b \in B$ with a neighborhood of $p \in X_b$.

This being true for general $p \in
\X$, it follows that $Z = B \times \psi_b(X_b)$ everywhere. Finally, since
the map $\psi$ is
assumed birational, it follows that $\X$ is  the normalization of $Z$; thus
it is likewise a
product, the map $\psi = id_B \times \psi_b$ and the Kodaira-Spencer map
identically zero.
\qed

We  now introduce the map

$$
\overline{\k _b}:T_bB \la H^0(X_b,\N_b/\Nt )
$$
defined to be the composition of $\k _b$ with the natural map $H^0(\N_b
)\to  H^0(X_b,\N_b/\Nt )$.
Notice that the Lemma implies that such a  map is an injection on the
subspace $\Im (\k _b)$.

As an application, we will fix the argument given above for plane curves.
In fact, we will prove  a slightly more general result. Let $Y$ be a smooth
surface and $D$ an effective  divisor on $Y$. Let
$V$ be an irreducible component of the Severi variety
of curves of given geometric genus $g$  that are linearly equivalent to $D$.
There is a universal family $\U \subset V\times Y$ of curves
over $V$, consider the normalization $\U ^{\nu}$ of $\U$ and
let $B\subset V$ be an open
subset over which $\U ^{\nu }$ is smooth.
Let $\X \to B$  be the restriction of $\U ^{\nu} $ to $B$ and
$\psi :
\X
\to B
\times Y$ be  as usual.
Then we obtain the following well known  result (cf. for example [K] and
loc. cit.) .

\proclaim Corollary 2.4. If, for general $b \in B$ we have
$
\deg _W(\psi_b^*\omega_Y) \; < \; 0
$
on every component $W$ of $X_b$, then
$$
\dim B \; \le \; -\deg(\psi_b^*\omega_Y) + g - 1 \, .
$$
In particular,  $\dim V_{d, \delta } = 3d + g -1$.

\ni {\it Proof}. We have
$$
\dim B \; \le \; \dim\left(\Im(\kappa_b)\right)  \le h^0(X_b,
\N_b/(\N_b)_{{\rm tors}})  \, $$
where the last inequality follows from the previous Lemma.

By the  definition of $\N _b$, we obtain
$$
\deg(c_1(\N_b)) \; = \; -\deg(\psi_b^*\omega_Y) + \deg(\omega_{X_b})
$$
hence, by our hypothesis
$\deg(c_1(\N_b)) >\deg(\omega_{X_b})
$ on each component of $X_b$. We now state
for both present and future use the following simple corollary of the
Riemann-Roch theorem for
curves:

\proclaim Observation 2.5. Let $X$ be a smooth curve of genus $g$, and $L$
any line bundle on
$X$ of degree $d$ such that
$L \otimes \omega_X^{-1}$ has positive degree on each component of $X$. If
$M$ is any line
bundle
 on $X$ such that
$L \otimes M^{-1}$ has nonnegative degree on each component of $X$, then
$$
h^0(X, M) \; \le \; h^0(X, L) \; = \; d-g+1 \, .
$$
If moreover the line bundle $L \otimes \omega_X^{-1}$ has degree  2 or more
on each
component of $X$, then $h^0(X, M)=d-g+1$ only if $\deg M= \deg L$.

\

Applying this to  $X = X_b$ and the line bundles $L = \psi_b^*\omega_Y^{-1}
\otimes \omega_{X_b}$ and $M = \N_b/(\N_b)_{{\rm tors}}$, we have
$$
\eqalign{\dim B \; &\le \; h^0(X_b,
\N_b/(\N_b)_{{\rm tors}}) \cr
&\le \; h^0(X_b, \psi_b^*\omega_Y^{-1}
\otimes \omega_{X_b}) \cr
&= \; -\deg(\psi_b^*\omega_Y) + p_a(X_b) - 1 \, . \cr}
$$
\qed

\

We  need now  to consider deformations of a map $X_b \to Y$ that
preserve tangency conditions with a fixed smooth curve $G \subset Y$. There
are two cases,
depending on whether we require tangency at a fixed point $p \in G$ or
allow tangency at a
variable point.

So,  let  $Y$ be a smooth surface, $G \subset Y$ a curve and $p \in G$ a
smooth point. Let
$\X
\to B$ be as above a smooth family of curves over a reduced base
$B$, $\psi : \X \to B \times Y$ a morphism of $B$-schemes, and $Q \subset
\X$ a section such
that  the pullback divisor
$$
\psi^*(G) - mQ \; \ge \; 0 \, .
$$
Let $b \in B$ be a general point and $q = X_b \cap Q$; suppose $\psi(q)=p$.
Let $v \in T_bB$  and
let
$\s =
\k_b(v)
\in H^0(X_b, \N_b)$ the corresponding first-order deformation, and $\ol \s
=\ol {\k _b}$.
Suppose finally that the differential
$d\psi_b$ vanishes to order
$l-1$ at
$q$, so that the image $\psi_b(\D)$ of a small neighborhood $\D$ of $q \in
X_b$ will have
multiplicity $l$ at $p = \psi(q)$. We have then the

\

\proclaim Lemma 2.6.  Let $\ol \s \in \Im \ol{\k _b}$.
Then
 $\ol \s$ vanishes to order at least $m-l$ at $q$, and cannot vanish to
order exactly $k$
for any $k$ with
$m-l<k<m$.
Moreover,  if we assume that
$\psi (Q)$ is a point,  $\ol \s$ vanishes to order at least $m$ at $q$.

\

\ni {\it Proof}. It will be sufficient to do this in case $B$ is
one-dimensional. Next, since $B$ is reduced and $b \in B$ is general, we
may assume $B$
smooth at $b$; so that, restricting to an analytic  neighborhood of $b \in
B$ we may
take $b$ the origin in an open subset $B$ of the affine line $\A^1 = \Spec
k[\e]$. Finally, from the
statement of the Lemma it is enough to prove it in case the divisor
$\psi^*(G)$ contains the curve
$Q$ with multiplicity exactly $m$; and since again $b \in B$ is general we
may assume as well
that the divisor $\psi_b^*G$ on $X_b$ contains the point $q$ with
multiplicity exactly $m$
as well.

Now, choose coordinates $(x,y)$ in an analytic  neighborhood of $p =
\psi(q)$ so that the
curve $G$ is given simply as the zero locus of $y$.
Let then   ${\partial \over \partial x}$ and ${\partial \over
\partial y}$ be the generators of the rank $2$ bundle $TY$ at $p$;
 we will abuse notation  and write ${\partial \over \partial x}$ and
${\partial \over
\partial y}$ also for the corresponding sections of $\psi _b^*T_Y$.

The first thing we will show is that  the image  of ${\partial \over
\partial x}$ in
$\N_b/(\N_b)_{{\rm tors}}$  vanishes to order $m-l$ at $q$.

We treat the case $l <m$ first for simplicity, and leave the case $l=m$ for
later.
Let $t$ be an $m^\th$
root of $\psi_b^*y$ in a neighborhood of $q \in \X_b$, then $t$ will be a
local coordinate on
$X_b$ near $q$ and the map $\psi_b$ will be given as
$$
\psi_b \; : \; t \; \longmapsto \; (t^l + c_{l+1}t^{l+1} + \ldots,\,t^m)
$$
so that the differential $d\psi_b$ is  given by
$$
\eqalign{d\psi_b \; : \; {\partial \over \partial t} \; &\mapsto \;
(lt^{l-1} + (l+1)c_{l+1}t^l +
\ldots){\partial \over
\partial x} + mt^{m-1}{\partial \over \partial y} \cr
&= \; t^{l-1}\left((l + (l+1)c_{l+1}t +
\ldots){\partial \over
\partial x} + mt^{m-l}{\partial \over \partial y}\right) \, , \cr}
$$
Denote
$
\tau(t)  : =  (l + (l+1)c_{l+1}t +
\ldots){\partial \over
\partial x} + mt^{m-l}{\partial \over \partial y}
$;
we have that the torsion subsheaf $(\N_b)_{{\rm tors}} \subset \N_b$ is
isomorphic to
$\O_{X_b}/\gm_q^{l-1}$, generated by the section $\tau (t)$.
Moreover, the quotient
$$
\N_b/(\N_b)_{{\rm tors}} \; = \; \O_{X_b}\{{\partial \over
\partial x}, {\partial \over \partial y}\}/\langle \tau \rangle
$$
 is generated for example by the image of the section ${\partial
\over \partial y}$.  Note finally that modulo  the subsheaf generated by $\tau$,
$$
{\partial
\over \partial x} \; \sim \; {mt^{m-l} \over l + (l+1)c_{l+1}t +
\ldots}\cdot {\partial
\over \partial y}
$$
so that the image of the
section ${\partial \over
\partial x}$ in $\N_b/(\N_b)_{{\rm tors}}$  vanishes to order exactly $m-l$
at $q$.

Now,  a general deformation $\psi$ of the map $\psi_b$ over the base $B
\subset \A^1_\e$ may be given in terms of coordinates $t$ and $\e$ on $\X$
near $q$ as

$$
\psi (t,\e)  = \bigl(\e ; \;  t^l  + c_{l+1}t^{l+1} +
\ldots + \e(\a_0+\a_1t+\ldots) + (\e)^2,
 t^m + \e(\b_0+\b_1t+\ldots) + (\e)^2 \bigr) \, .
$$

The condition that the divisor $\psi^*((y)) = mQ$ near $q$ says that
we can take $t$ to be an $m^\th$ root of the pullback $\psi^*y$ not just on
$X_b$, but in a
neighborhood of $q$ in $\X$. This means that a deformation satisfying the
hypotheses of the
lemma may be written as
$$
\psi (t,\e)  = \bigl( \e ,t^l  + c_{l+1}t^{l+1} +
\ldots + \e(\a_0+\a_1t+\ldots) + (\e)^2,
 t^m \bigr) \, .
$$

From the definitions, the image $\k_b({\partial
\over \partial \e}) \in H^0(X_b, \N_b)$ of the tangent vector ${\partial
\over \partial \e} \in T_bB$ under the Kodaira-Spencer map  will
be given as the image in $\N_b$ of
$$
\s : = \; \k_b({\partial
\over \partial \e}) \; = \; (\a_0+\a_1t+\ldots){\partial \over
\partial x} \, ,
$$
whose image $\ol \s$ in $\N_b/(\N_b)_{{\rm tors}}$, as we have seen,
vanishes to order at least
$m-l$ at
$q$. Moreover, since $b \in B$ is general, the differential $d\psi_\e$ will
vanish to order $l-1$ at
$X_\e \cap Q$ for all $\e$ near $b$; that is, $t^{l-1}|d\psi_\e$. This
implies that
$$
\a_1 = \a_2 = \ldots = \a_{l-1} = 0 \, ;
$$
or in other words, $\ol \s$
cannot vanish to order exactly
$m-l+1,\ldots,m-1$ at $p$. To complete the proof in case $m > l$, the
further condition that
$\psi(Q)
\equiv p
\in G$ says that
$\a_0=0$, which further implies that $\ol \s$ vanishes to order at least
$m$ at $q$.

The case $m=l$ is  completely analogous.  As before we write the map $\psi_b$ as
$$
\psi_b \; : \; t \; \longmapsto \; (t^n + c_{n+1}t^{n+1} + \ldots,\,t^m)
$$
where now $n \geq m$.  We leave it to the reader to check that the same
argument yields that
if $\psi(Q) \equiv p \in G$,  the section $\ol \s$ vanishes to order at
least $m$ at
$q$.
\qed

\

\ni  2.3. {\bf Dimension counts and consequences.}
We will now use the general theory developed above to establish
Propositions 2.1 and
2.2.

\ni { \it Proof of Proposition 2.1}. To begin with, it follows from the
naive dimension
count at the beginning  that the Severi variety
$\V$ has dimension at least $2d+g-1+|\b|$.

We thus have to show that $\dim \V \le 2d+g-1+|\b|$ everywhere. To do this,
let $V \subset
\V$ be an irreducible component, $\X \subset V \times \P^2$
the universal family curve over $V$, and $\X^\nu$ the normalization of the
total space. We will
actually
restrict our attention to the open subset of $ V$ over which $\X^\nu$ is
smooth,  which we will still call $V$.

Let $[X] \in V$ be a general point, so that the restriction $\nu =
\psi_{[X]}$ of $\psi$ to
the fiber of
$\X^\nu$ over
$[X]$ is the normalization $\nu : X^\nu
\to X
\subset
\P^2$  of the corresponding curve $X \subset \P^2$; and let $\N$ be the
normal sheaf of the map $\nu$;
notice that this might appear as an abuse of notation, as we have already
used the symbol $\N$ with a different
meaning, in the previous chapter; we hope that this will not create
confusion.  By the definition of
$\V$, we have in an analytic neighborhood of $[X]$ a collection of $|\a|$
and $|\b|$ sections
$\{Q_{i,j}\}$ and
$\{R_{i,j}\}
\subset \X^\nu$ such that
$$
\psi(Q_{i,j}) \; = \; p_{i,j}
$$ and
$$
\psi^*(L) \; = \; \sum i\cdot Q_{i,j} \; + \; \sum  i\cdot R_{i,j} \, .
$$
Let $q_{i,j} = Q_{i,j} \cap X^\nu$ and $r_{i,j} = R_{i,j} \cap X^\nu$.
Note that the points $\{q_{i,j}\}$ are necessarily distinct, since they
have distinct images
$p_{i,j} \in L \subset \P^2$. We may assume as well that the points
$\{r_{i,j}\}$ are
distinct, and disjoint from the $\{q_{i,j}\}$: if not, $[X]$ being general
in $V$, $V$ would be as
well a component of a Severi variety $V^{d,\d}(\a',\b')$ for some
$(\a',\b')$ with $|\b'| < |\b|$,
which we will show has dimension $2d+g-1+|\b'| < 2d+g-1+|\b|$.

We need to introduce one more bit of notation. We
denote by $l_{i,j}-1$ the order of vanishing of the differential $d\nu$ at
the point
$r_{i,j}$. We then let $D$ and $D_0 \in
\Div(X^\nu)$ to be the divisors
$$
D \; = \; \sum_{1 \le j \le \a_i} i \cdot q_{i,j} + \sum_{1 \le j \le \b_i}
(i-1) \cdot r_{i,j} \, .
$$
and
$$
D_0 \; = \; \sum_{1 \le j \le \a_i} (l_{i,j}-1) \cdot r_{i,j}
\, .
$$
Note that $D$ is a divisor of degree
$$
\eqalign{\deg(D) \; &= \; I\a + I\b - |\b| \cr
&= \; d-|\b| \, . \cr}
$$
and that
$$
\deg\left((\nu^*\O_{\P^2}(1))(-D)\right) \; \ge \; 0 \, .
$$
Note also that
$
\deg\left(c_1(\N_{{\rm tors}})\right) \; \ge \; \deg(D_0) \, ,
$
on every component of $X^\nu$, with equality holding if and only if $\nu$ is an
immersion away from $\{r_{i,j}\}$; so that
$$
\deg\left(c_1(\N/\N_{{\rm tors}})\right) \; \le \;  \deg(c_1(\N)) -
\deg(D_0)
$$
again with equality holding if and only if $\nu$ is an
immersion away from $\{r_{i,j}\}$.

Finally, let $D_1$ be the effective part of $D - D_0$.

Now, applying Lemmas 2.3 and 2.6, we see that
$$
\dim \V \; \le \; h^0(X^\nu, (\N/\N_{{\rm tors}})(-D_1)) \, .
$$
We have
$$
\deg\left((\N/\N_{{\rm tors}})(-D_1)\right) \; \le \;  \deg(c_1(\N)) -
\deg(D)
$$
and since
$$
c_1(\N) \; = \; \nu^*\O_{\P^2}(3) \otimes \omega_{X^\nu}
$$
we see that the line bundle
$$
(c_1(\N)(-D)) \otimes
\omega_{X^\nu}^{-1}  \; = \; \left((\nu^*\O_{\P^2}(1))(-D)\right) \otimes
\nu^*\O_{\P^2}(2)
$$
has strictly positive degree on each component of $X^\nu$. We may thus
apply the simple
Observation 2.5 to the line bundles $c_1(\N)(-D)$ and $(\N/\N_{{\rm
tors}})(-D_1)$ to conclude that
$$
\eqalign{\dim \V \; &\le \; h^0(X^\nu, (\N/\N_{{\rm tors}})(-D_1)) \cr
&\le \; \deg\left(c_1(\N)(-D)\right) - g + 1 \cr
&= \; (3d+2g-2-\deg(D)) - g + 1 \cr
&= \; 2d+g-1+|\b| \, . \cr}
$$
\qed

\ni \underbar{Remark}.  Notice that the argument above implies that the
image  of the Kodaira Spencer
map can be identified as follows:
$$
\Im \ol{ \k }_{[X]} = H^0(X^\nu, (\N/\N_{{\rm tors}})(-D_1))
$$

\

\ni {\it Proof of Proposition 2.2}. We  start by establishing what is
perhaps the
subtlest point: that the map
$\nu$ is indeed an immersion. In fact, much of this has already been
accomplished in the proof
of Proposition 2.1. Keeping the notations introduced there, we
see that since the line bundle $(c_1(\N)(-D)\otimes
\omega_{X^\nu}^{-1}$ on $X^\nu$ has degree at least 2 on any component of
$X^\nu$,
we may apply the  second part of
Observation 2.5 to deduce the equality
$$
(\N/\N_{{\rm tors}})(-D_1) \; = \;  c_1(\N)(-D)
$$
so that $D_1 = D - D_0$ and
$$
\N/\N_{{\rm tors}} \; = \;  c_1(\N)(-D_0)
$$
and hence $\nu$ is an
immersion away from $\{r_{i,j}\}$.

To see that $\nu$ is an immersion at the point $r_{i,j}$, we may assume
that the component
$X_0$ of
$X^\nu$ containing $r_{i,j}$ does not map to a line, so that the line
bundle $(c_1(\N)(-D)\otimes
\omega_{X^\nu}^{-1}$ has degree at least 4 on $X_0$. It follows that there
exists a section
$\ol \s$ of
$c_1(\N)(-D) = (\N/\N_{{\rm tors}})(-D_1)$ vanishing to order exactly 1 at
$r_{i,j}$; and
by the previos Remark, this section  must be in the image of the
Kodaira-Spencer map
$$
\k_{[X]} \; : \; T_{[X]}V \; \la \; H^0(X^\nu, (\N/\N_{{\rm tors}})).
$$
 But the multiplicity of $r_{i,j}$ in the divisor $D_1 = D - D_0$ is
$(i-1)-(l_{i,j}-1) = i -
l_{i,j}$, and it follows that $\ol \s$, viewed as a section of $\N / \Nt $,
 vanishes to order exactly $i -
l_{i,j}+1$ at $r_{i,j}$. By Lemma 2.6, then, we must have $l_{i,j} = 1$;
that is, $\nu$ must be an
immersion at $r_{i,j}$

To show that $X$ has only
nodes as singularities, we have to show it has no triple points and that no
two branches are
tangent to each other. For the former, if $s, t, u \in X^\nu$ are points
mapping to the same
point $p \in X \subset \P^2$, it is enough to show that there exists a
section of
$\N(-D)$ vanishing at $s$ and $t$ but not at $u$. This follows immediately from
Riemann-Roch: if $s$ and $t$ and $u$ all belong to the same component of
$X^\nu$, that
component must map to a plane curve of degree at least 4, so that $\N \otimes
\omega_{X^\nu}^{-1}$ will have degree at least $8$ there; while if two lie
on the same
component, $\N(-D) \otimes \omega_{X^\nu}^{-1}$ will have degree at least
$6$ there.
Similarly, for the latter, it is enough to show that  if $s, t \in X^\nu$
are points mapping to
the same point $p \in X \subset \P^2$,  there exists a section of the sheaf
$\N(-D)$ vanishing
at $s$ but not at $t$, which follows from the same argument.

\

As for parts $c$ and $d$ of 2.2, we have already seen just
from the dimension statement that the points $\{q_{i,j}\}$ and $\{r_{i,j}\}$ are
all distinct, since otherwise $V$ would be a component of a Severi variety
$V^{d,\d}(\a',\b')$
for some $(\a',\b')$ with $|\b'| < |\b|$; and the same logic implies that
the points $\{s_{i,j}\}$
are disjoint from the points $\{p_{i,j}\}$. To see that the points
$s_{i,j}$ are all distinct, on the
other hand, it is sufficient to observe that, by the argument of the
preceding paragraph, for
any
$(i',j')
\ne (i,j)$ with $s_{i',j'} = s_{i,j}$, there is a section of the sheaf
$\N(-D)$ vanishing at $r_{i,j}$ but not at $r_{i',j'}$; by Lemma 2.6 this
will correspond to a
deformation of $X$ in which $s_{i',j'}$ moves but $s_{i,j}$ stays still.

\

Next, given that $\nu : X^\nu \to X$ is an immersion,  part
$b$ follows from $d$; if
 $\nu$ is one-to-one over points of $L$ then $X$ is smooth along $L$.

\

Finally, part $e$: if a branch of $X$ corresponding to a point $s \in
X^\nu$ were tangent to $G$, it would be enough to show that there exists a
section of $\N(-D)$
vanishing at $s$, which we know; and likewise if two points $s, t \in
X^\nu$ mapped to the
same point $p \in G$, it would be enough to show that there exists a
section of $\N(-D)$
vanishing at $s$ but not at $t$, which again we know.
\qed

\

The following  restatement of Proposition 2.1 will be useful in the
applications in
the next section. To set it up, fix a line $L \subset \P^2$ and a finite
subset $\Omega \subset L$.
Let $V \subset \P^N$ be any irreducible, locally closed subset of the space
of plane curves of
degree
$d$, and $[X] \in V$ a general point. Let
$\pi :
W
\to
X \subset \P^2$ be any map not constant on any irreducible component of
$W$, whose degree
over each irreducible component $X_i$ of $X$ is equal to the multiplicity
of $X_i$ in $X$ (so
that in particular the pullback
$\pi^*\O_{\P^2}(1)$ has degree
$d$). Let $g$ be the geometric genus of $W$, and let $e$ be the cardinality
of the intersection
$\#\left(X \cap (L \setminus \Omega)\right)$. We have then

\

\proclaim Corollary 2.7.
$$
\dim V \; \le \; 2d+g-1+e \, ;
$$
and if
equality holds and $\#\left(X \cap (L \setminus \Omega)\right) =
\#\pi^{-1}(L \setminus
\Omega)$ then
$V$ is a dense open subset of a generalized Severi variety $\V$.

\

\ni {\it Proof}. This follows readily from (2.1), after a few reductions.
To begin with,
it is enough to prove this in case $W$ is smooth, since replacing $W$ by
its normalization only
strengthens the inequality. Secondly, it is enough to do it in case $X$ is
irreducible: applying
the statement to the inverse image of each component of $X$ in turn and
adding the results
yields the desired inequality in general. (This second reduction is not
really essential, but
will allow us to refer to the degree of the map $W \to X$ without confusion.)

Now, since $W$ is smooth, the map $W \to X$ factors through the
normalization $X^\nu \to
X_{\rm red}$, and we claim that it is enough to prove it in case $W =
X^\nu$. To see this,
assume the result proved in case $W = X^\nu$ and  consider what happens if
$W \to X^\nu$ is a
finite map of degree $m > 1$. In this case the degree of $X_{\rm red}$ is
$d/m$, and the
genus $h$ of $X^\nu$ is related to the genus $g$ of $W$ by Riemann-Hurwitz:
$$
g \; \ge \; mh - m + 1
$$
Now, applying (2.1) directly to $X_{\rm red}$, we have
$$
\eqalign{\dim(V) \; &\le \; 2{d \over m} + h - 1 + e \cr
&\le  \; 2{d \over m} + {g-1 \over m} + e \cr
&< 2d+g-1+e  \cr}
$$
and we have a contradiction.

We may thus assume that $W = X^\nu$. Now, suppose first that $e=0$. In this
case,  the
statement we want to prove is exactly Proposition 2.1 and we are done. More
generally,
consider the map $\phi$ from a neighborhood of $[X] \in V$ to $|\O_L(e)|$
sending a point $[W]
\in V$ to the reduced intersection $W \cap (L \setminus \Omega)$. Applying
the $e=0$ case of
the statement  of the Corollary to the fiber of $\phi$ over a general point
$D \in |\O_L(e)|$
(replacing $\Omega$ by $\Omega' =
\Omega \cup {\rm supp}(D)$), we conclude that the fibers of $\phi$ have
dimension at most
$2d+g-1$, and hence that $\dim(V) \le 2d+g-1+e$.
\qed

Note that, in the case of equality, the map $\pi : W \to X$ is necessarily
a birational isomorphism
on each component of $W$.

\

\ni  2.4. {\bf Normal sheaves and normal bundles.}
To conclude this section, we should say a few words about the relationship
between the
 treatment of Severi varieties given here and other possible approaches.
Briefly,
there are two ways of analyzing the deformations of a plane curve $X$ satisfying
certain geometric conditions. In the approach taken here, which we may call
the ``parametric"
approach, we look at deformations of the normalization map $\nu: X^\nu \to
X \subset \P^2$;
so that the tangent space to the space of deformations is a priori a
subspace of the space of
sections of the normal sheaf $\N$ of the map. This has the virtue (at
least, it is a virtue in our
present circumstances) of incorporating the condition that the geometric
genus of $X$ is
preserved in the deformations. Moreover, the sheaf $\N$ is a sheaf on a
smooth curve. On the
other hand, it has the defect that, until we know that
$\nu$ is an immersion, the sheaf
$\N$ may have torsion.

In the other approach, which we will call the ``Cartesian" approach, we
look instead at
deformations of
$X$ as a subscheme of $\P^2$;
so that the tangent space to the space of deformations is a priori a
subspace of the space of
sections of the normal bundle $N_{X/\P^2} \cong \O_X(d)$ of the divisor $X
\subset \P^2$.
This is in some ways more direct---all we are doing, after all, is
practicing the time-honored
tradition of varying the coefficients of the defining polynomial of
$X$---and it is in
particular useful when we want to intersect our family with other
subvarieties of the space
$\P^N$ of plane curves of degree $d$. But it has the drawback that we have
to impose extra
conditions to ensure that the geometric genus of
$X$ stays constant. These conditions, moreover, sometimes interact badly
with conditions such
as tangency with a fixed curve.

What is the relationship between the two? In case $\nu : X^\nu \to X$ is an
immersion, it is
reasonably straightforward. To start with, let $\I \subset \O_X$ be the
{\it conductor ideal} of
$X$. This may be characterized in several equivalent ways:
\medskip
\item{$\bullet$} \thinspace It is the annihilator of the sheaf
$\nu_*\O_{X^\nu}/\O_X$;
\smallskip
\item{$\bullet$}  \thinspace It is the largest ideal $\I \subset \O_X$ such
that the pullback
map $\nu^*$ gives a bijection between ideals in $\O_X$ contained in $\I$
and ideals in
$\O_{X^\nu}$ contained in $\nu^*\I$;
\smallskip
\item{$\bullet$}  \thinspace On an affine open subset  of $X$ with defining
equation $f(x,y)$, it
is the ideal of polynomials $g(x,y)$ such that the 1-form
$$
 \nu^*({g(x,y) dx \over \partial f / \partial y})
$$
is regular on $X^\nu$; and
\smallskip
\item{$\bullet$}  \thinspace More concretely, in case $\nu : X^\nu \to X$
is an immersion, it is
the ideal in $\O_X$ whose restriction to each branch $\Delta_i$ of $X$ at
each point $p \in X$
is equal to the restriction to that branch of the ideal of the union of all
other branches of $X$
through
$p$. In other words, if $p_i \in X^\nu$ is the point lying over $p$ in the
branch $\Delta_i$,
$$
\nu^*\I \; = \; \O_{X^\nu}\Bigl(- \sum_i \bigl(\sum_{j \ne i} {\rm
{mult}}_p(\Delta_i \cdot
\Delta_j)\bigr)\cdot p_i\Bigr) \, .
$$

However we characterize the conductor, it is not hard to see that,  in case
$\nu : X^\nu \to X$ is
an immersion, the normal sheaf $\N$ of the map $\nu$ and the normal bundle
$N_{X/\P^2}
\cong \O_X(d)$ of the curve $X$ are related by
$$
\N \; = \; \nu^*(\I \otimes N_{X/\P^2}) \, .
$$
This is perhaps most easily seen in terms of the last description of the
conductor: if the local
defining equation $f(x,y)$ of $X$ at a point $p \in X$ factors in the
completion of the local ring
$\O_{X,p}$ as
$$
f(x,y) \; = \; f_1(x,y)f_2(x,y)\cdots f_n(x,y)
$$
then a general first-order deformation of the map will simply move each
branch, resulting in a
curve given by the equation
$$
f_\e(x,y) \; = \; (f_1(x,y)+\a_1\e)(f_2(x,y)+\a_2\e)\cdots
(f_n(x,y)+\a_n\e) \, .
$$
As a deformation of the map, that is, as a section of  $\N$, this will be
nonzero at the point of $X^\nu$ corresponding over the branch $\Delta_i$
given by $f_i(x,y)=0$
if and only if the coefficient $\a_i \ne 0$. But the corresponding section
of the normal bundle,
that is, the restriction to $X$ of the coefficient of $\e$ in
$f_\e(x,y)$, on this branch is $\a_i\prod_{j \ne i}f_j(x,y)$, which
vanishes to order $\sum_{j \ne i}  {\rm {mult}}_p(\Delta_j
\cdot
\Delta_i)$.

In any event, the conclusion is that the sections of the  $N_{X/\P^2}$ coming
from deformations of the map are simply those lying in the conductor ideal
(or, classically,
``satisfying the adjoint conditions", in view of the third characterization
above). Moreover, if
we impose further conditions of tangency with fixed curves, the allowed
deformations of the
map correspond to sections of $\N$ vanishing to the appropriate order at the
points of $X^\nu$ lying over  the points of tangency; and these sections,
by the second
characterization above, correspond to sections of ${\cal J} \otimes
N_{X/\P^2}$ for a unique
ideal sheaf ${\cal J} \subset \I$.

We thus have a very useful dictionary between the two languages, at least
as long as $\nu$ is
an immersion. Otherwise the correspondence is more complicated. For
example, if $[X]$ is a point on the variety of plane curves of given degree
$d$ and genus
$g$, corresponding to a curve with a cusp and $\d-1={d-1 \choose 2} - g -1$
nodes, then in a
neighborhood of $[X]$ we cannot simultaneously normalize the fibers of the
universal family
$\X \subset V \times \P^2 \to V$; so deformations of $X$ preserving the
geometric genus do
not correspond  to deformations of the map.
\

\vfill\eject

\

\

\ni {\sc 3. Hyperplane sections of Severi varieties: set-theoretic description}.

\

We are now prepared to describe the hyperplane sections
of the generalized Severi varieties.
 In this chapter we will prove Theorem 1.2, showing that the intersection $\V
\cap H_p$ is indeed a union of generalized Severi varieties of dimension
one less, and
saying which ones potentially occur. In the following chapter we will prove
Theorem
1.3, establishing that all the generalized Severi varieties listed as possible
components of $\V
\cap H_p$ do in fact occur, and describing the multiplicities with which
they appear.

\

\ni 3.1. {\bf The basic setup.}
Let $V'$ be any
irreducible component of the intersection $\V
\cap H_p$ and $[X_0] \in V'$ a general point. If $X_0$ does not contain $L$
it is
easy to see that $V'$ must be a component of one of the generalized Severi
varieties
listed in the first part of Theorem 1.2, so we will focus on the case
$L
\subset X_0$. We then consider a curve
$\Gamma = \{[X_\g]\} \subset \V$ passing through the point $[X_0]$, and the
corresponding family of plane curves $\X \to \Gamma$. Applying a variant of
semistable reduction to this family in a neighborhood of $[X_0] \in \Gamma$
we arrive
at a family $\Y \to B$ of nodal curves dominating the curves $X_\g$ in our
family.
Analyzing this family, we find a number of  geometric conditions that the curve
$X_0$ must satisfy, which limit the number of its degrees of freedom.
 Playing these off against the fact that $X_0$ is a
general member of a  variety of dimension
$ \dim(\V)-1$ we are led to our
characterization of
$\V
\cap H_p$.

Although the approach is quite simple, the arguments  tend to appear
extremely complicated.  In fact, there are a priori no restrictions on  the
number of components, for example, of the special fiber  of the family $Y \to B$
or their configuration, the notation alone can be very cumbersome. We will
therefore
proceed in two steps: we will give the analysis first subject to a number of
simplifying assumptions, which will make the logic of the argument
relatively clear
(all of these hypotheses, moreover, will in fact turn out to satisfied in
reality). Then
we will go back and prove the result without assumptions. Comparing the argument
here with the first should make it clear why in fact the assumptions hold:
if any of
them were indeed violated, we could replace one of the inequalities of the first
calculation with a strict inequality, and so arrive at a contradiction.

We start by defining the families $\X \to \Gamma$ and $\Y \to B$ that we
will be working with, and establishing the relevant notation. As above, we
let $[X_0]
\in V'$ be a general point of a component $V'$ of $\V
\cap H_p$, and $\Gamma \subset \V$ a curve containing $[X_0]$. We will assume
that  the general point $[X_\g]$ of
$\Gamma$ is a general point of $\V$; that is, it satisfies the conclusions of
Proposition 2.1  above.

Now, let $\nu : \Gamma^\nu \to \Gamma$ be the normalization of $\Gamma$, and
choose a point $b_0 \in \Gamma^\nu$ lying over $[X_0]$. Let $\X^\nu$ be the
normalization of the total space of the pullback $X \times_\Gamma
\Gamma^\nu$, so
that the family $\X^\nu \to \Gamma^\nu$ has as general fiber a smooth curve of
genus $g = {d-1 \choose 2} - \d$.

Next, we want to carry out a nodal reduction of the family $\X^\nu \to
\Gamma^\nu$ in a neighborhood of $b_0$ to arrive at a family of nodal curves $\Y
\to B$
satisfying  the requirements below. (This can be achieved
 after possibly further base change and blowing
up of the nodal reduction.)
Let $\Y^* = f^{-1}(B \setminus
\{b_0\})$ be the complement of the special fiber $Y_0 = f^{-1}(b_0)$ of $\Y
\to B$.
Then we require that

\

\item{a.} The total space $\Y$ is smooth.
\smallskip
\item{b.} The map carrying a general fiber $Y_b$ of $\Y \to B$ to the
corresponding
plane curve $X_\g \subset \P^2$ extends to a regular morphism
$$
\pi \; : \; \Y \; \la \; \P^2 \, .
$$
\smallskip
\item{c.} The inverse image $\pi^{-1}(L) \cap \Y^*$ consists of $|\a| +
|\b|$ disjoint
sections
$\{Q^*_{i,j}\}_{1 \le j \le \a_i}$ and $\{R^*_{i,j}\}_{1 \le j \le \b_i}$,
with $\pi(Q^*_{i,j})
\equiv p_{i,j}$ and $R^*_{i,j}$ intersecting the general fiber in a point
$r_{i,j}$ of
multiplicity
$i$ in the divisor $(\pi|_{Y_b})^*L$---that is we have
$$
\pi^*L \cap \Y^* \; = \; \sum i \cdot Q^*_{i,j} + \sum i \cdot R^*_{i,j} \, .
$$
\smallskip
\item{d.} The closures $Q_{i,j}$
and
$R_{i,j}$ in $\Y$ of the sections $Q^*_{i,j}$
and
$R^*_{i,j}$ are still disjoint
and  they do not pass through
any of the singularities of $Y_0$.
\smallskip
\item{e.} Finally, $\Y \to B$ is minimal with respect to these properties.

In sum, we have the following diagram of objects and morphisms:

\

\vskip2.2in

\hskip.5in \special {picture semistablediagram}

\

\ni 3.2. {\bf Some simplifying assumptions and some corollaries.}
In order to present as clearly as possible the actual picture of the
families $\X \to \Gamma$ and
$\Y \to B$,  we will first carry out the analysis of the family $\Y \to B$ under
three simplifying assumptions, all of which we will show in the last part
of this section do
in fact hold. We will also mention  some interesting  facts that
will  follow
as consequences of the proof of Theorem 1.2 (in particular they will not be
assumed in the course
of the proof).Since they contribute to the picture of the family $\Y \to B$
we will state
them here.

The first of our assumptions is perhaps the least obvious: it is that

\

\ni
{\bf Assumption (a).} {\it The curve $X_0$ contains $L$ with multiplicity
1; that is, $X_0 = X \cup L$,
where $X$ is a plane curve of degree $d-1$ not containing $L$.}

\

\ni (We should remark that this seems to be
false in slightly more general situations, for example if we consider the
variety
$V$ of plane curves of given degree and geometric genus having a triple point.)

Now, given (a), we see that the special fiber $Y_0 = f^{-1}(b_0)$ of $\Y \to B$ will
contain a unique component $\tilde L$ such that $\pi$ maps $\tilde L$ onto
$L$, and
indeed the map $\pi|_{\tilde L} : \tilde L \to L$ will be an isomorphism.
We may then
group together the remaining components of $Y_0$ into two sets: we will let $Y
\subset Y_0$ be the union of the irreducible components mapping to $X$ on
which $\pi$ is
nonconstant, and $Z \subset Y_0$ the union of the irreducible components of
$Y_0$ on which
$\pi$ is constant.  In these terms, we will assume next that

\

\ni
{\bf Assumption (b).} {\it The curve $Z \subset Y_0$ consists of a disjoint
union of chains of rational
curves joining $\tilde L$ to $Y$.}

\

\ni
{\bf Assumption (c).}  {\it The sections $\{Q_{i,j}\}_{1 \le j \le \a_i}$
and $\{R_{i,j}\}_{1 \le j \le
\b_i}$ are disjoint from $Z$.}

 \

Note that by the last statement, each such section meets either $Y$ or
$\tilde L$
but not both. To keep track of how many of the sections $Q_{i,j}$ pass
through each, we will introduce some more notation. First, we define two
further sequences
$\a'$ and
$\a''$ with
$\a' +
\a'' =
\a$: we let
$\a'_i$ be the number of the sections $\{Q_{i,j}\}_{j=1,\ldots,\a_i}$
passing through
$Y$, and $\a''_i$ the number passing through $\tilde L$. We will likewise
label the
sections passing through $Y$ (respectively, $\tilde L$) as $\{Q'_{i,j}\}_{1
\le j \le \a'_i}$
(respectively,
$\{Q''_{i,j}\}_{1 \le j \le \a''_i}$), their points of intersection with
$Y_0$ as $\{q'_{i,j}\}_{1 \le j
\le \a'_i}$ (respectively,
$\{q''_{i,j}\}_{1 \le j \le \a''_i}$), and their image points in
$L$ as the subset
$\Omega' =
\{p'_{i,j}\}_{1 \le j \le \a'_i}$ (respectively, $\{p''_{i,j}\}_{1 \le j
\le \a''_i}$).

\

As for the sections
$R_{i,j}$, the situation is a little different, by virtue of the first of
our corollaries:

\proclaim  Consequence 1.
Every
section $R_{i,j}$ passes through $Y$.

Notice that this means that $\b \leq \b '$, so that this statement
 is actually part of  Theorem 1.2.
Thus, we do not introduce a  new
set of symbols. Rather, we will provisionally let $\b^0_i$ be the number of
the sections
$\{R_{i,j}\}_{j=1,\ldots,\b_i}$ passing through
$Y$, and suppose (after possibly relabeling) that the sections $R_{i,j}$
passing through $Y$ are
 $\{R_{i,j}\}_{1 \le j \le \b^0_i}$. We will let $r_{i,j}$ be the point of
intersection of
$R_{i,j}$ with $Y_0$, and $s_{i,j} = \pi(r_{i,j}) \in L$ its image point.

Note also that, given (b), we may, at the expense of introducing rational
double points into our
surface
$\Y$, collapse the connected components of $Z$ to points, so that the fiber
$Y_0$ consists
simply of
$\tilde L$ and $Y$ (we have done this in the diagram below for clarity). We
will index
the points of intersection of $\tilde L$ with $Y$ as follows: for each $i$,
we let
$r''_{i,1},\ldots,r''_{1,\b''_i}$ be the points of $\tilde L \cap Y$
appearing with multiplicity $i$ in
the divisor
$\pi^*L|_Y$. We will also let $s''_{i,j} = \pi(r''_{i,j}) \subset L$ be the
images of these points.

In sum, we have the following picture of the family $\Y \to B$ and its sections:

\vfill\eject

\

\vskip3.4in

\special {picture semistable}

\

\ni (Note that we have anticipated in this picture the statement (d) that every
section $R_{i,j}$ passes through $Y$.)

We also want to underline two other corollaries of the proof:

\

\proclaim Consequence 2.
The curve $X$ is reduced; that is, the map $\pi|_Y : Y \to X$ is a birational
isomorphism on each component of $Y$; and

\proclaim Consequence 3.
 $Y$ is smooth (though not in general  connected).

\

\ni 3.3. {\bf Proof of  simplified Theorem 1.2}
We now  prove  Theorem 1.2 under  assumptions (a), (b) and (c).
We do that by comparing two  different relations that
the arithmetic  genus $g$ of $Y$  has to satisfy

To begin with, note that the curves $Y$ and $\tilde L$ intersect in
$|\b''|$ points, which
are nodes of $Y_0$. Since the arithmetic genus of $\tilde L \cong L \cong
\P^1$ is 0, we have
the first relation
$$
g \; = \; p_a(Y_0) \; = \; p_a(Y) + |\b''| - 1 \, .
$$

Now we apply Corollary
2.7 to obtain a second relation. We see that the dimension of the family in
which $X$ can move is at most $2(d-1) + g(Y) -
1 + |\b'|$. But $X$ is a general member of the
$V'$, which has dimension
$$
\dim(V') \; = \; \dim(\V)-1 \; = \; 2d+g-2+|\b| \, .
$$
Thus
$$
\eqalign{2d+g-2+|\b|  \; &\le \; 2(d-1) + g(Y) - 1 + |\b'| \cr
&\le \; 2(d-1) + p_a(Y) - 1 + |\b'| \cr
&= \; 2d-2+g-|\b''| + |\b'|  \cr
&= \; 2d-2+g  + |\b^0| \cr
&\le \; 2d+g-2+|\b| \, .\cr}
$$
We conclude that equality holds throughout, and that $V'$ is therefore a
component of the
Severi variety $V^{d-1,\d'}(\a', \b')(\Omega')$ satisying the equality
$$
g' \; = \; g(X) \; = \; g - |\b' - \b| + 1
$$
or equivalently
$$
\eqalign{\d - \d' \; &= \; ({d-1 \choose 2} - g) - ({d-2 \choose 2} - g') \cr
&= \; d-2 + g - g' + 1 \cr
&= \; d-1 - |\b' - \b| \, . \cr}
$$

Another way to view this situation is via the
pullback
$\pi^*(L)$ of $L$ to $\Y$, and in particular its restriction to $Y$.
We have
$$
\pi^*L \; = \; m\cdot \tilde L + D + \sum_{1 \le j \le \a_i} i \cdot
Q_{i,j} + \sum_{1 \le j \le \b_i} i
\cdot R_{i,j}
$$
 where $D$ is supported on $Z$ and $m$ is some integer. Restricting to $Y$
we have
$$
\pi^*L|_Y \; = \; \sum_{1 \le j \le \a'_i} i \cdot q'_{i,j} + \sum_{1 \le j
\le \b^0_i} i \cdot r_{i,j} +
\sum_{1 \le j \le \b''_i} i \cdot r''_{i,j} \, .
$$
In particular, if we set $\b' = \b^0 + \b''$, the degree
$$
d-1 \; = \; \deg(\pi^*L|_Y) \; = \; I\a' + I\b' \, .
$$
The point is, of the $d-1$ points of intersection of $X$ with $L$, $I\a'$
will occur at the
assigned points $p'_{i,j}$, while the remaining $I\b' = I\b^0 + I\b''$ will
occur at the
$|\b'|$ unassigned points $s_{i,j}$ and $s''_{i,j}$. The greatest degree of
freedom would thus
seem to be attained when $\a'$ is as small as possible and $\b' = \b^0 +
\b''$ as large as possible.
But
$\b^0$ is bounded above by $\b$, and there is a penalty for taking $\b''$
large: this will decrease
the geometric genus of the curve
$X$, which will drop the dimension of the family which it can move.

\

Note some additional consequences of this analysis. First, we
see from the fact that equality holds in the last of the series of
inequalities above that
$\b^0 =
\b$, so that  every section $R_{i,j}$ passes through $Y$, as stated, in
other words, $\beta \le \beta '$.
Statements (e) and (f)
above likewise follow from the proof: the fact that $X$ is reduced is a
consquence of the
application of Corollary 2.7; while if $Y$ were singular, we would have
$g(Y) < p_a(Y)$, giving rise
to a strict inequality in the second of series of inequalities above.
This completes the proof of Theorem 1.2 subject
to hypotheses (a)-(c).

\

\ni 3.4. {\bf The local picture of the degeneration.}
Based on the above analysis, we have a complete picture of the behavior of
the family of plane
curves $X_\g$ as they degenerate to $X_0$. Away from $L$, there is no apparent
degeneration; the family is equisingular. We will describe the family near
each of the
relevant points of $L$:

\

\item{$\bullet$} \thinspace At a point $p''_{i,j}$---that is, a point of
$\Omega \setminus
\Omega'$---the curve $X_0$ is smooth, since $X$ does not pass through
$p''_{i,j}$.
Thus, in a neighborhood of $p''_{i,j}$ we see the curves
$X_\g$ simply ``flatten out" to the line $L$.

\

\item{$\bullet$} \thinspace At a point $p'_{i,j} \in \Omega'$, the
curve $X_0$ has an $i^\th$ order tacnode, since $X$ is smooth at $p'_{i,j}$
and has
contact of order $i$ with $L$ there. On the other hand, the inverse image of
$p'_{i,j}$ in $Y_0$ has two distinct points, one in $Y$ and one in $\tilde
L$; in other
words, the map $Y_0 \to X_0$ factors through the normalization of $X_0$ at
$p_{i,j}$. Thus, in an analytic neighborhood $U$ of $p'_{i,j}$ the curves
$X_\g$ will have
two branches, one tending to a neighborhood of $p'_{i,j}$ in $L$ and one to a
neighborhood of $p'_{i,j}$ in $X$. Moreover, these two branches of $X_\g$
will have $i$
points of intersection in $U$, merging to form the one point of
intersection multiplicity
$i$; thus, as the curves $X_\g$ approach $X_0$ we see $i$ nodes of the curves
$X_\g$ approach the point $p'_{i,j}$ and coalesce to form the tacnode of
$X_0$. Note finally
that, of the two branches of $X_\g$ near $p'_{i,j}$, it is the one tending to a
neighborhood of $p'_{i,j}$ in $X$ that has contact of order $i$ with $L$ at
$p'_{i,j}$.

\vskip1.3in

\special {picture p'degeneration}

\item{$\bullet$} \thinspace The picture at a point $s_{i,j} \in L$---that
is, a limit of an
unassigned point of intersection of $X_\g$ with $L$---is exactly the same
as the picture
above near a point $p'_{i,j}$: the nearby curves $X_\g$ will have two
branches in a
neighborhood of $s_{i,j}$, one tending to $L$ and one to $X$; and
correspondingly we
will see $i$ nodes of the curve $X_\g$ merge to form the $i^\th$-order
tacnode of
$X_0$ at $s_{i,j}$. The only difference, in fact, is that where in the
preceding case the
curves $X_\g$ all had a fixed point of contact of order $i$ with $L$ at
$p'_{i,j}$, in this
case the curves $X_g$ have a point of intersection multiplicity $i$ with
$L$ at a
variable point tending to $s_{i,j}$. Note that as in the preceding case, of
the two branches of
$X_\g$ near $s_{i,j}$, it is the one tending to  $X$ that has contact
of order $i$ with $L$ at $s_{i,j}$.

\

\item{$\bullet$} \thinspace The picture near the ``new" tacnodes of
$X_0$---that is,
the points $s''_{i,j}$ of intersection of $X$ with $L$ that are not limits
of points of
intersection of $X_\g$ with $L$---is the most interesting. Here, the
inverse image of
$s''_{i,j}$ in $Y_0$ is a single point, which is an ordinary node of $Y_0$.
It follows that
in an analytic neighborhood of $s''_{i,j}$ the curves $X_\g$ are
irreducible, with $i-1$
nodes tending to the $i^\th$-order tacnode $s_{i,j}$ of $X_0$. Local
equations for (and another
picture of) such a family will be given in section 4.1 below, when we
discuss the
geometry of deformations of a tacnode.

\vskip1.4in

\special {picture s''degeneration}

Note that we can now account for all the nodes of the general curve $X_\g$
of our
family as $X_\g$ tends to $X_0$. First, $\d'$ of the nodes stay away from
$L$, and
become simply nodes of $X \subset X_0$. At each point $p'_{i,j}$, we see
$i$ nodes of
$X_\g$ absorbed into the tacnode formed by $X$ and $L$; and similarly at
$s_{i,j}$.
Finally, we see $i-1$ nodes absorbed into each point $s''_{i,j}$. In
particular, we see
once more that the number $\d'$ of nodes of $X$ is
$$
\d' \; = \; \d - I\a' - I\b - (I\b'' - |\b''|)
=  \; \d - I\a' - I\b' + |\b''|
= \; \d - (d-1) + |\b' - \b| \, .
$$

We can also use the above analysis to give a complete description of the
components of the
curve $Z$. Let us suppose that the curve $\tilde L$ appears with
multiplicity $m$ in the
pullback divisor $\pi^*L$.  Since the
restriction of $\pi^*L$ to $Y$ has multiplicity $i$ at each point
$r''_{i,j}$, if $i < m$ it cannot meet
$\tilde L$ itself: rather, it must meet a component $Z_1$ of $Z$ that
appears with multiplicity
$i$ in
$\pi^*L$. Since the degree of the restriction $\pi^*L|_{Z_1}$ is zero and
$Z_1$ has
self-intersection $-2$, we see that
$Z_1$ must meet another component of $Y_0$ that appears with multiplicity
$2i$ in $\pi^*L$,
and so on. Ultimately, we see that $i$ must divide $m$, and that $Z$ will
include a chain of
$m/i - 1$ rational curves joining $r''_{i,j}$ to  $\tilde L$, which appear
with multiplicities $i, 2i,
3i,\ldots,m-2i, m-i$ in $\pi^*L$.

\

The picture of the actual fiber $Y_0$ of $\Y \to B$ is thus:

\vfill\eject

\

\vskip 2.6in

\hskip .5in \special {picture Y0actual}

\

\ni where the length of the chain of rational curves joining $\tilde L$ to
a point $r''_{i,j} \in Y$ is
$m/i - 1$. Note that  the multiplicity $m$ with which curve $\tilde L$
appears in the
pullback divisor $\pi^*L$ must therefore be a multiple of $\lcm(\b'')$. We
will see in the following
section that for suitably general $\Gamma$, we have $m = \lcm(\b'')$.

\

\ni 3.5. {\bf Verifying the assumptions.}
To complete the proof of Theorem 1.2, we will now go back and  justify
the assumptions
(a)-(c). This  amounts to setting up the same calculation without any
of the hypotheses, and then observing that if any of them were  violated we
would have a
strict inequality in one of the series of inequalities above.

To begin with, we  extend the definitions of $\tilde L$, $|\b^0|$, $|\b''|$
and $|\b'|$ to the
general setting. First, we take
$\tilde L$ to be simply the union of the components of the fiber $Y_0$
dominating $L$. Then, we
define $\b^0_i$ to be the number of the sections
$\{R_{i,j}\}_{j=1,\ldots,\b_i}$ that meet either
$Y$ or a connected component of $Z$ meeting $Y$---again, if we consider the
limits as $\g \to 0$ of
the $|\b|$ unassigned points of intersection of $X_\g$ with $L$, $|\b^0|$
will be simply the number
 that lie on $X$. Similarly, we let
$|\b''|$ be the number of points of $Y$ meeting $\tilde L$, plus the number
of connected
components of $Z$ meeting both $\tilde L$ and $Y$, and set $|\b'| = |\b^0|
+ |\b''|$.

A key observation is that no connected component of $\pi^{-1}(L)$ can be
contained in $Z$, since otherwise we could contract it to obtain an
isolated point in the
inverse image of $L$. In other words, every connected component of $Z$
mapping to a point of
$L$ must meet $\tilde L$. Thus, with these definitions, we have as before
$$\#\left(X
\cap (L
\setminus
\Omega)\right) \; \le \; |\b'| \, .
$$
Moreover, the genus satisfies
$$
g = p_a(Y_0) \; \ge \; p_a(\tilde L) + p_a(Y) + |\b''| - 1.
$$

Now consider the hypothesis (a) that $X_0$ contains $L$ simply.
If this were not the case, that is,  $X_0 = m\cdot L \cup X$, where $X$ is now a
plane curve of degree
$d-m$. Since the map
$\pi|{\tilde L} : \tilde L \to L$ has degree $e$, we have to replace the
equality
$p_a(\tilde L) = 0$ with the inequality
$$
p_a(\tilde L) \; \ge \; -m+1
$$
and correspondingly in place of $p_a(Y) = g - |\b''| + 1$ we have
$$
p_a(Y) \; \le \; g - |\b''| + m \, .
$$
This appears to work against us. But at the same time, the degree of the image
$X$ of $Y$ is now $d-m$ rather than $d-1$, so that the dimension of the
family in which
it moves is
$$
\eqalign{2(d-m) + p_a(Y) + |\b'| - 1 \; &\le \; 2(d-m) + g - |\b''| + |\b'|
+ m - 1 \cr
&= \;  2d+g+|\b| -m - 1 \cr
&< \; 2d+g+|\b| - 2  \cr}
$$
but our hypothesis is that $X$ moves in $V'$ that has dimension $2d+g+|\b|
- 2$, hence
 we have a contradiction.

\

Assumption (b) is more intuitively clear, if slightly more cumbersome to
check. The point is, our
basic inequality on the genus of $Y$ was based on the fact that
$$
g \; = \; p_a(Y_0) \; = \; p_a(Y) + p_a(\tilde L) + |\b''| - 1 \, .
$$
Now, if any connected component of $Z$ met $L$ twice, or met $Y$ twice, or
itself had
strictly positive arithmetic genus, we would have a strict inequality
$$
g \; > \; p_a(Y) + p_a(\tilde L) + |\b''| - 1 \,
$$
and so would arrive at a contradiction. Thus each connected component $Z_0$
of $Z$ is a tree of
rational curves meeting each of $\tilde L$ and $Y$ at most once.

To verify (b), then, we simply
have to check that every leaf of this tree (that is, every irreducible
component $W$ of $Z$ meeting
at most one other component of $Z$) meets either
$Y$ or
$\tilde L$. But, if it did not, by the minimality of $\Y$ at least two of
the sections $Q_{i,j}$ and
$R_{i,j}$ would have to meet it; otherwise we could blow down $W$ in $\Y$
and still satisfy all the
conditions imposed on $\Y$. Now, clearly no two of the sections $Q_{i,j}$
can meet the same
component of $Z$, since they have distinct images in $\P^2$. On the other
hand, our basic
estimate on the dimension of $V'$ is based on the fact that the curve $X =
\pi(Y)$ can have at
most $|\b'| = |\b| + |\b''|$ points of intersection with $L$ outside of
$\Omega$, corresponding to the
points of intersection of $Y$ with the sections $R_{i,j}$ and the points of
intersection of $Y$ with
the connected components of $Z$ meeting $\tilde L$. If a section
$Q_{i,j}$ and a section $R_{i,j}$ met the same component of $Z$, the point
$s_{i,j} = \pi(r_{i,j})$
would lie in $\Omega$, and there would be strictly fewer than $|\b'|$
points of intersection of
$X$ with $L \setminus \Omega$, and by Corollary (2.7) the dimension of $V'$ would be strictly
less than $2d+g+|\b| - 2$.  Likewise, if two of the sections
$R_{i,j}$ met the same component of
$Z$, two of the points
$s_{i,j}$ would coincide, and
$X$ would again have strictly fewer than $|\b'|$ points of intersection
with $L$ outside of
$\Omega$.

\

Finally, the verification of (c) follows the same pattern as that of (b):
given that $Z$ consists
simply of chains joining $\tilde L$ to the points $r_{i,j} \in Y$, if any
of the sections
$Q_{i,j}$ met $Z$ at all, we would have $s_{i,j} \in \Omega$; while if any
of the sections $R_{i,j}$
met $Z$ at all, we would have $s_{i,j} = s''_{i', j'}$ for some $i,j,i',j'$
and once more $X$ would
again have strictly fewer than $|\b'|$ points of intersection with $L$
outside of
$\Omega$.

\vfill\eject

\

\

\ni \cl{{\sc 4. Hyperplane sections of Severi varieties: local geometry}}

\

In this section we  describe the geometry of  $\V$ in a
neighborhood of a general point of a component of its intersection with a
hyperplane
$H_p$.  As might be expected, this is relatively straightforward for
$V^{d,\d}(\a+e_k,\b-e_k)$ (the first possibility described in Theorem 1.2),
and substantially more complex in the
case of a component of  $V^{d-1,\d'}(\a',\b')$. In the
analysis of the latter case, we will rely heavily on the description given
in [CH] of
various loci in the versal deformation space of an $m$th order tacnode. The
arguments in this section are based on conversations with Ravi Vakil, and
appear as
well in [V].

\

\ni  4.1. {\bf Deformation spaces of tacnodes.}
We start by recalling the relevant results from [CH]
(the reader is referred to Section 2.4 of [CH] for  details). To begin
with, let $(C,p)$ be an
$m^\th$ order tacnode, that is, a curve singularity analytically equivalent
to the
origin in the plane curve given by the equation
$$
y(y-x^m) \; = \; 0 .
$$
The versal deformation space of $(C,p)$ is then the family $\pi : \S \to
\Delta$, where
$\Delta \cong \A^{2m-1}$ with coordinates
$(a_{m-2},\ldots,a_0,b_{m-1},\ldots,b_0)$,
$\S$ is the subscheme of $\Delta \times \A^2$ given by the equation
$$
\eqalign{f(&x,y,a_{m-2},\ldots,a_0,b_{m-1},\ldots,b_0) =\cr
&= \; y^2 + (x^m +
a_{m-2}x^{m-2}+ \cdots +a_1x + a_0)y + b_{m-1}x^{m-1} + \cdots + b_1x + b_0 \cr
&= \; 0
\cr}
$$
and $\pi : \S \to \Delta$ is the projection $\S \subset \Delta \times \A^2
\to \Delta$.

As is clear from the description of the generalized Severi varieties in the
preceding
section, we are primarily interested in two loci in the base $\Delta$ of
the versal
deformation: the locus $\Delta_m \subset \Delta$ of points $(a,b) \in
\Delta$ over
which the fiber $\S_{a,b}$ of $\S \to \Delta$ is reducible---equivalently,
the closure of
the locus of $(a,b)$ such that $\S_{a,b}$ has $m$ nodes---and  the closure
$\Delta_{m-1}$ of the locus of $(a,b)$ such that
$\S_{a,b}$ has $m-1$ nodes. Since the defining equation $f(x,y,a,b)$ for
$\S$ exhibits
$\S$ as a double cover of the $x$-line $\Delta \times \A^1_x$ over
$\Delta$, we can
describe these two loci in terms of the branch divisor of the cover:  the
discriminant
$\d =
\d_{a,b}(x)$ of the equation
$f(x,y,a,b)$ above as a quadratic polynomial in
$y$ is given by
$$
\d_{a,b}(x) \; = \; (x^m +
a_{m-2}x^{m-2}+ \cdots +a_1x + a_0)^2 - 4(b_{m-1}x^{m-1} + \cdots + b_1x +
b_0) ,
$$
and the loci $\Dm$ and $\Dn \subset \D$ are  the closure of the loci of $(a,b)$
such that $\d_{a,b}$ has $m$ double roots and such that $\d_{a,b}$
has $m-1$ double roots, respectively ($\Dm$ could also be characterized  as the
locus of squares). In particular, we see that
$\Dm$ is simply the
$(m-1)$-plane $\Dm \cong \A^{m-1} \subset \D$ given by the equations $b_{m-1} =
\cdots = b_1 = b_0 = 0$; and that $\Dn$ is an $m$-dimensional subvariety of $\D$
containing
$\Dm$ (and having multiplicity $m$ at a general point of $\Dm$).

In these terms, we can now state the results of [CH] that we will use here:

\proclaim Lemma 4.1. Let $m \ge 2$, and let $W \subset \D$ be any smooth,
$m$-dimensional subvariety containing the $(m-1)$-plane $\Dm$, and suppose only that
its tangent plane at the origin is not contained in the hyperplane $H
\subset \D$
given by $b_0 = 0$. Then we have
$$
W \cap \Dn \; = \; \Dm \cup \Gamma
$$
where $\Gamma$ is a smooth curve having contact of order exactly $m$ with $\Dm$
 at the origin.

There is an alternative way to express this lemma, which is very useful in
both its
proof and the present appliction; it may also provide more geometric
insight. Let $\tilde
\D$ be the blow-up of
$\D$ along the
$(m-1)$-plane
$\Dm$; let
$E \subset \tilde \D$ be the exceptional divisor  and $F \cong \P^{m-1}
\subset E$ the
fiber of
$E$ over the origin in $\D$. Let $\tilde \Dn$ be the proper transform of
$\Dn$ in $\tilde
\D$. Then we can state the last result as the

\proclaim Lemma 4.2. The intersection $\tilde \Dn \cap E$ contains $F$ as a
component of
multiplicity $m$. Moreover, $\Dn$ is smooth at any point of $F$ not
contained in the
proper transform $\tilde H$ in $\tilde \D$ of the hyperplane $H \subset \D$
given by $b_0 = 0$.

To see the equivalence of the two statements, note that by the second the
tangent
plane to
$\Dn$ at any point of
$F$ not contained in the proper transform of the hyperplane $b_0=0$ must contain
the tangent plane to $F$ and be contained in the tangent plane to $E$. In
particular,
the proper transform $\tilde W$ of any subvariety $W \subset \D$ satisfying the
hypotheses of the first statement must intersect $\Dn$ transversely in a smooth
curve $\tilde \Gamma$ having intersection multiplicity $m$ with $E$ at its
point of
intersection with $F$; and the first statement follows. Conversely, if we
apply the first
statement just to the $m$-planes $W_b$ in $\D$ containing $\Dm$, we see
that $\Dn$
contains $F$ and (away from the proper transform of the hyperplane $b_0=0$) is
swept out by the smooth curves $\tilde W_b$; the second statement follows.

\vfill\eject

\

\vskip3.2in

\hskip.3in \special {picture tildeDn}

\

\

It will be helpful also to recall the steps in the proof of these lemmas in
[CH]. Briefly,
the statement of Lemma 4.1 is first verified by direct calculation in case
$W$ is the
specific
$m$-plane
$b_{m-1} = b_{m-2} = \cdots = b_1=0$. Then the homogeneity of $F \setminus
(F \cap
\tilde H)$ under the action of the automorphism group of the singularity $(C,p)$
allows us to deduce it for $W$ any $m$-plane containing $\Dm$ and not
contained in
$H$. This, as noted, is all we need to deduce the statement of Lemma 4.2,
and then
Lemma 4.1 in general follows as well.

\

\ni 4.2. {\bf Products of deformation spaces of tacnodes.}
What we have to do now is to use this information to develop an analogous
picture in
the product of deformation spaces associated to a collection of
higher-order tacnodes.
For what follows, then, we will let $m_1, m_2, \ldots$ be any sequence of
integers
$m_j \ge 2$, and $(C_j,p_j)$ be an
$(m_j)^\th$ order tacnode. We will denote the versal deformation of
$(C_j,p_j)$ by
$\pi_j : \S_j
\to
\Delta_j$, and let
$(a_{j, m_j-2},\ldots,a_{j,0},b_{j,m_j-1},\ldots,b_{j,0})$ be coordinates
on $\Delta_j$ as
above. For each $j$, we will let $\Delta_{j,m_j}$ and $\Delta_{j, m_j-1}
\subset \Delta_j$
be as above the closures of the loci in $\Delta_j$ over which the fibers of
$\pi_j$ have
$m_j$ and
$m_j-1$ nodes respectively. Finally, we set
$$
\Delta \; = \; \Delta_1 \times \Delta_2 \times \cdots ,
$$
$$
\Delta_m \; = \; \Delta_{1,m_1} \times \Delta_{2,m_2} \times \cdots ,
$$
and
$$
\Delta_{m-1} \; = \; \Delta_{1,m_1-1} \times \Delta_{2,m_2-1} \times \cdots .
$$
Note that $\Delta$, $\Dm$ and $\Dn$ have dimensions $\sum 2m_j-1$, $\sum
m_j-1$ and
$\sum m_j$ respectively.

Our goal now is to describe how a smooth subvariety $W \subset \Delta$ of
dimension
$\sum (m_j-1) + 1$, containing $\Dm$, will intersect $\Dn$, again with some
hypothesis
on its tangent plane at the origin. Specifically, let $H \subset \D$ be the
union of the
hyperplanes $(b_{j,0} = 0) \subset \D$, and suppose that the tangent plane
to $W$ is not
contained in $H$. By the dimension count, we would expect
$W$ to intersect $\Dn$ in the union of $\Dm$ and a residual curve $\Gamma$; we will
show that this is indeed the case, and that the intersection number of
$\Gamma$ with
$\Dm$ at the origin in $W$ is $\prod m_j$. (What will be different from the
single-tacnode case is the local geometry of $\Gamma$: as we will see, it
may have many
branches, each of which may be singular at the origin.) To make the full
statement,
let $\l$ the the least common multiple of the $m_j$, let $\mu =
\prod m_j$ and set
$\k =
\mu / \l$. We will prove the

\proclaim Lemma 4.3. With the hypotheses above, in an
\'etale neighborbood of the origin in
$\D$  the intersection
$$
W \cap \Dn \; = \; \Dm \cup \Gamma_1  \cup \Gamma_2 \cup \ldots \cup \Gamma_\k
$$
where $\Gamma_1,\ldots,\Gamma_\k \subset W$ are distinct reduced unibranch
curves having intersection multiplicity exactly
$\l$ with
$\Dm$
 at the origin. The curves $\Gamma_\a$ will all be smooth if $\l = m_j$ for
some $j$;
otherwise they will all be singular, with multiplicity $\l/\max_j\{m_j\}$.

As before, it will be helpful to express this in terms of the geometry of a
blow-up: we
let $\tilde \Delta$ be the blow-up of $\Delta$ along the plane $\Dm$, $E \subset
\tilde \Delta$ the exceptional divisor, $F
\subset E$ the fiber of
$E$ over the origin in $\D$. Let $\tilde \Dn$ be the proper transform of
$\Dn$ in
$\tilde
\D$. Then we can state the last result as the

\proclaim Lemma 4.4. The intersection $\tilde \Dn \cap E$ contains $F$ as a
component of
multiplicity $m$. Moreover,  in an \'etale neighborhood of any point $p \in
F$ not
contained in the proper transform $\tilde H$ of $H$ in $\tilde \D$, $\tilde \Dn$
consists of $\k$ reduced branches, each having multiplicity $\l/\max_j\{m_j\}$,
intersection number $\l$ with $E$ along $F$, and tangent cone at $p$
supported on a
linear space contained in $E$.

\ni {\it Proof}. We will prove this in two stages: first, we
will verify the  statement of Lemma 4.3 directly in case $W$ is any linear
subspace of
$\D$; from this we will deduce  the full statement of Lemma 4.4, and hence
Lemma 4.3 for arbitrary
$W$. Note that we can not do this as in the single-tacnode case by giving
an explicit
calculation in the case of a particularly simple plane $W$ and then invoking
homogeneity, for one reason: in the single-tacnode case we have an apriori
lower bound
of $m$ for the intersection multiplicity of $\Gamma$ with $\Dm$, so that if
we verify
that $\Gamma$ is a smooth curve with contact of order $m$ with $\Dm$ for one
plane $W$, we may deduce it for any plane $W'$ such that $W$ lies in the
closure of the
orbit of $[W'] \in G(m,2m-1)$ under the action of the automorphism group of the
deformation. Here we do not have the analogous lower bound $(\Dm \cdot
\Gamma) \ge
\mu$, and so we have to deal first with an arbitrary linear space $W$.

So: let $W \subset \D$ be any plane of dimension $\sum (m_j-1) + 1$ containing
$\Dm$ and not contained in $H$;  $W$ will be spanned by $\Dm$
and one additional vector $v \in \D$. Let $t$ be a nonzero linear function
on $W$
vanishing on the hyperplane $\Dm$. Let
$\rho_j :
\D
\to
\D_j$ be the projection, and
$v_j = \rho_j(v)$, so that the image $W_j = \rho_j(W) \subset \D_j$ of $W$
in $\D_j$ will
be the plane spanned by the subspace $\D_{j,m_j}$ and the vector $v_j$. Since by
hypothesis
$v_j$ does not lie in the hyperplane $b_{j,0} = 0$, by Lemma ** of [CH], we
may write the
intersection
$W_j
\cap
\D_{j,m_j}$ as the union of
$\D_{j,m_j}$ and a smooth curve $\Gamma_j$ having contact of order $m_j$ with
$\D_{j,m_j}$ at the origin. It follows that in some \'etale neighborhood of
the origin in
$W_j$ we may choose coordinates $(x_{j,0},x_{j,1},\ldots,x_{j,m-2},t_j)$ so
that
$$
(\rho_j)^*t_j \; = \; t ;
$$
and the hyperplane
$\D_{j,m_j}$ is given by $t_j=0$ and the curve $\Gamma_j$ by the equations
$$
x_{j,1} \; = \; x_{j,2} \; = \; \cdots \; = \; x_{j,m-2} \; = \; 0
$$
and
$$
(t_j)^{m_j} \; = \; x_{j,0} .
$$
We may then take the collection of functions
$\{y_{j,i} = (\rho_j)^*x_{j,i}\}_{1 \le i \le m_j-2}$ and $t$ as local
coordinates in an \'etale
neighborhood of the origin in $W$, in terms of which the the hyperplane $\Dm$ is
given by $t=0$ and the residual intersection $\Gamma$ of $W$ with $\Dn$ by the
equations
$$
y_{j,i} \; = \; 0 \; ,  \qquad \forall \; j \; {\rm and} \; i \; : \; 1 \le
i \le m_j-2
$$
and
$$
t^{m_j} \; = \; y_{j,0}  \; ,  \qquad \forall \; j \; .
$$

$\Gamma$ is a curve, since specifying the value of the coordinate $t$ at a
point of
$\Gamma$ determines  the value of the coordinates $y_{j,0}$ (and hence all the
coordinates $y_{j,i}$) up to a choice of an $(m_j)^\th$ root. More
explicitly, suppose that
we choose for each $j$ an $(m_j)^\th$ root $\zeta_j$. Then we may
parametrize a branch $\Gamma_\zeta$ of the curve $\Gamma$ by
$$
t \; = \; z^\l ,
$$
$$
y_{j,0} \; = \; {z^{\l/m_j} \over \zeta_j}
$$
and of course $y_{j,i}=0 \; \forall \; i>0$. This parametrization is
one-to-one, since the
powers of
$z$ appearing have no common factor; the multiplicity of the image at the
origin is the
smallest power of $z$ appearing, which is $\l/\max\{m_j\}$; and since the
pullback of $t$
is $z^\l$, the intersection multiplicity of the image with the hyperplane
$\Dm \subset W$
defined by $t=0$ is $\l$. Moreover, all branches of $\Gamma$ are
parametrized in this
fashion; and two collections of roots
$\zeta_j$ and
$\eta_j$ will give rise to the same branch if, and only if, for some
$\l^\th$ root of unity
$\e$ we have
$$
\zeta_j \; =  \; \e^{m_j}\eta_j
$$
for all $j$. Since $\zeta_j = \e^{m_j}\zeta_j$ for all $j$ only if $\e =
1$, the number of such
branches is the number $\mu$  of collections of roots $\zeta_j$ divided by
$\l$, that is,
$\k$. Thus the statement of Lemma 4.3 is established for any linear space $W$.

The remainder of the argument for Lemmas 4.3 and 4.4 is straightforward:
exactly as
before, the statement of Lemma 4.3 for  a linear space $W$ (satisfying the
hypotheses of
the Lemma) implies the statement of Lemma 4.4, which in turn implies 4.3 in
general. To
carry this out, let $U = \P(\D/\Dm) \setminus \P(H)$ be the complement of the
projectivizations of the hyperplanes $b_{j,0} = 0$ in the projectivization
of the
quotient $\D/\Dm$, so that we have a morphism
$$
\tau \; : \; V \; = \; \tilde \D \setminus \tilde H \; \la \; U
$$
expressing the complement $V$ of the proper transform $\tilde H$ of $H$ in
the blow up
$\tilde \D$ as a projective bundle (with fiber dimension $\sum m_j -1$)
over $U$.
Let
$\tilde
\Dn^0 =
\tilde
\Dn
\cap V$ be the intersection of the proper transform $\tilde \Dn$ with this
open subset of
$\tilde
\D$, and let $\sigma : \tilde \Dn^0 \to U$ be the restriction of $\tau$ to
$\tilde \Dn^0$.
The statement of Lemma 4.3 for linear spaces
$W$ says precisely that the fibers of
$\s$ are curves consisting of $\k$ reduced branches, each having multiplicity
$\l/\max_j\{m_j\}$, intersection number $\l$ with $E$ at its unique point $p$ of
intersection with $F$, and tangent cone at
$p$ supported on a linear space contained in the tangent space to $E$. It
follows that
exactly the same is true of $\tilde \Dn$ in a neighborhood of any point $p
\in F$ not in
$\tilde H$: it has
$\k$ reduced branches, each having multiplicity
$\l/\max_j\{m_j\}$ and intersection multiplicity $\l$ with $E$ along $F$,
and tangent
cone  supported on a linear space contained in the tangent space to $E$. In
other words,
we have proved Lemma 4.4; and as before Lemma 4.3 follows for an arbitrary
smooth
$m$-dimensional subvariety $W \subset \D$ satisfying the hypotheses of the
Lemma.
\qed

\

\ni 4.3. {\bf The local geometry around irreducible curves.}
In the remaining two parts of this section we will complete the proof of
Theorem 1.3
by analyzing the geometry of a generalized Severi variety $\V$ in a
neighborhood of a
general point $[X_0]$ of one of the generalized Severi varieties $V'$
listed in Theorem
1.2. We will start in this subsection with the (relatively) simple case of
a general point
 of
$V^{d,\d}(\a+e_k,\b-e_k)(\Omega \cup
\{p\})$; in the following one we will apply the preceding results to carry
out the
analysis at a general point of $V^{d-1,\d'}(\a',\b')(\Omega')$.

\

So: assume that $\b_k > 0$, and let $[X_0]$ be a general point of $V' =
V^{d,\d}(\a+e_k,\b-e_k)(\Omega
\cup
\{p\})$. We have then the

\

\proclaim Proposition 4.5. The variety $\V$ contains $[X_0]$; it is smooth
there  and has
intersection multiplicity
$k$ with $H_p$ along $V'$.

\ni {\it Proof}. This follows directly from an analogous statement about the
linear series of divisors on $L \cong \P^1$. To set this up, consider the
rational map
$$
\pi \; : \; |\O_{\P^2}(d)| \cong \P^N \; \la \; |\O_L(d)| \cong \P^d
$$
given by restriction (note that this is a linear projection, with vertex
the subspace in
$\P^N$ of curves containing $L$). Inside the target space $|\O_L(d)|$, we
consider three
loci: we let $H = \{D : D - p \ge 0\}$ be the hyperplane of divisors
containing the point
$p$; and we set
$$
\Phi \; = \; \left\{ D \in |\O_L(d)| \; : \; D = \sum_{1\le j \le \a_i}
i\cdot p_{i,j} +
\sum_{1\le j
\le \b_i} i\cdot p'_{i,j} \quad \hbox{for some  } \; p'_{i,j} \in L \right\}
$$
and similarly, setting $\b' = \b - e_k$,
$$
\Psi \; = \; \left\{ D \in |\O_L(d)| \; : \; D = k\cdot p + \sum_{1\le j
\le \a_i} i\cdot p_{i,j}
+
\sum_{1\le j
\le \b'_i} i\cdot p'_{i,j} \quad \hbox{for some  } \; p'_{i,j} \in L \right\} .
$$
Now, let $V = V^{d, \d } \subset \P^N$ be the ordinary Severi variety. we have
$$
H_p \; = \; \pi^{-1}(H) \; ;
$$
and if we let
$$
\s \; = \;  \pi|_V \; : \; V \; \la \; |\O_L(d)|
$$
be the restriction of $\pi$ to $V$,
$$
\V \; = \; \s^{-1}(\Phi)
$$
and
$$
V^{d,\d}(\a+e_k,\b-e_k) \; = \; \s^{-1}(\Psi) \; .
$$
Proposition 4.5 will thus follow from the combination of the two Lemmas

\proclaim Lemma 4.6.  The differential $d\s$ of $\s : V \to |\O_L(d)|$ is
surjective at $[X]$.

\ni and

\proclaim Lemma 4.7.   In a neighborhood of the point $D_0 = X \cdot L \in
|\O_L(d)|$,
the variety $\Phi$ is smooth and has intersection multiplicity $k$ with the
hyperplane $H$ along $\Psi$.

\ni {\it Proof of Lemma 4.6}. Since the map $\pi : |\O_{\P^2}(d)| \cong \P^N \to
|\O_L(d)| \cong \P^d$ is a linear projection, we need only check that the
projective
tangent plane $PT_{[X]}V \subset \P^N$ to
$V$ at
$[X]$ intersects the vertex $L + |\O_{\P^2}(d-1)| \cong \P^{N-d-1}$ of
$\pi$ transversely,
that is, in codimension $d+1$ in $PT_{[X]}V$. But now the projective
tangent space
$PT_{[X]}V$ is simply the linear series of curves of degree $d$ passing
through the $\d$
nodes of $X$, and its intersection with the vertex the linear series of curves
of degree $d-1$ passing through the nodes; and since the nodes impose
independent
conditions on curves of degree at least $d-2$ these will have dimensions
${d(d+3) \over
2} -
\d$ and ${(d-1)(d+2) \over
2} -
\d$ = ${d(d+3) \over
2} -
\d - (d+1)$ respectively. \qed

\ni {\it Proof of Lemma 4.7}. We can do this simply in coordinates: let $x$
be an
affine coordinate on $L \cong \P^1$ such that the point $p$ is given by
$x=0$; let the
point $p_{i,j}$ have coordinate $\l_{i,j}$ and suppose that the $\b'_i$
points other than
$p$ and $\{p_{i,j}\}$ at which the divisor $D_0$ has multiplicity $i$ have
coordinates
$\mu_{i,j}$ (note that by Lemma ** the $\l_{i,j}$ and the $\mu_{i,j}$ are
all distinct).
Then we can parametrize a neighborhood of $[D_0]$ in $\Phi$ by
$$
(\e, \e_{i,j}) \; \longmapsto \; [f(x)] = [(x-\e)^k \prod (x-\l_{i,j})^i \prod
(x-\mu_{i,j}-\e_{i,j})^i]
$$
from which we see in particular that $\Phi$ is smooth at the point $[D_0]$.
Now, writing a
point
$[f(x)]
\in |\O_L(d)|$ as
$$
f(x) \; = \; x^d + b_{d-1}x^{d-1} + \ldots + b_1x +
b_0
$$
the defining equation of the hyperplane $H \subset |\O_L(d)|$ is simply
$b_0 = 0$, which
pulls back via this parametrization to $\e^k$ times a polynomial in the
$\e_{i,j}$ nonzero
at the origin; it follows that in a neighborhood of $[D_0]$, the divisor
cut on $\Phi$ by
$H$ is simply $k$ times $\Psi$. \qed

\

\ni 4.4. {\bf The local geometry around reducible curves.}
It remains to describe the local geometry of a generalized
Severi variety $\V$ in a neighborhood of a point $[X_0]$ , where $X_0 = X
\cup L$ and
$[X]$ is a general point of a generalized Severi variety
$V' = V^{d-1,\d'}(\a',\b')(\Omega')$. Thus, we will suppose that $\Omega_i' =
\{p'_{i,j}\}_{1\le j \le \a'_i}$ is any subset of cardinality $|\a'|$ of
$\Omega_i =
\{p_{i,j}\}_{1\le j \le \a_i}$ such that $\{p'_{i,1},\ldots,p'_{i,\a'_i}\}
\subset
\{p_{i,1},\ldots,p_{i,\a_i}\}$ for each $i$; and that there are points
$\{q'_{i,j}\}_{1 \le j
\le \a'_i}$ and
$\{r'_{i,j}\}_{1 \le j \le \b'_i}$ in the normalization $\nu : \tilde X \to
X \subset \P^2$
such that $\nu(q'_{i,j}) = p'_{i,j}$ and the pullback
$\nu^*(L) =
\sum iq_{i,j} + \sum ir_{i,j}$. With this said, our basic result is the

\proclaim Proposition 4.8. In a neighborhood of $[X_0]$, the variety $\V$
will have
$$
{\b'
\choose
\b}{I^{\b'-\b} \over \lcm(\b' - \b)}
$$
branches, each of which will  have intersection multiplicity
$\lcm(\b'-\b)$ with $H_p$ along $V'$.

\ni {\it Proof}. As in the case of Lemma 4.5, we want to deduce this from a
local
calculation, in this case, Lemma 4.3. Before we can do this, we have to specify
which of the points $r_{i,j} \in \tilde X$ will be limits of points of
unassigned tangency
on nearby curves in the family; each such specification will determine a
collection of
branches of $\V$. So: to start with, choose any  subset $\Lambda =
\{r_{i,j}\}_{1\le j \le
\b_i}$ of the set $\{r'_{i,j}\}_{1\le j \le \b'_i}$ such that
$\{r_{i,1},\ldots,r_{i,\b_i}\}
\subset
\{r'_{i,1},\ldots,r'_{i,\b'_i}\}$ for each $i$. By way of notation, let
$\b'' = \b' - \b$, and
label the complement of the subset $\{r_{i,j}\} \subset \{r'_{i,j}\}$ as
$\{r''_{i,j}\}_{1 \le
j \le \b''_i}$; and let
$s_{i,j} \subset L \subset \P^2$ (respectively, $s'_{i,j}$, $s''_{i,j}$) be
the images of the
point $r_{i,j}$ (respectively, $r'_{i,j}$, $r''_{i,j}$). Similarly, set
$\a'' = \a - \a'$, and label
the complement of the subset $\{p'_{i,j}\} \subset \{p_{i,j}\}$ as
$\{p''_{i,j}\}_{1 \le j \le
\a''_i}$.

Now, in an analytic neighborhood of the point
$[X_0] = [X+L]
\in
\P^N$, we will define the {relaxed local Severi variety} $W_\Lambda$ to be the
closure of the locus of curves
$X_t$  satisfying the following six conditions:

\

\item{i)} $X_t$ preserves the $\d'$ nodes of $X$; that is, for every node
of $X_0$
away from $L$, $X_t$ will have a node nearby.
\smallskip
\item{ii)} At each point $p''_{i,j}$, $X_t$ has contact of order $i$ with $L$.
\smallskip
\item{iii)} In a neighborhood of each point $p'_{i,j}$, $X_t$ has $i$ nodes
\smallskip
\item{iv)} In a neighborhood of each point $s_{i,j}$, $X_t$ has $i$ nodes

\

To specify the remaining two conditions we need to make one remark.
Conditions iii)
 requires that, in an
\'etale or analytic  neighborhood of a point $p'_{i,j}$, the deformation
$X_t$ of $X$ will be
reducible; that is, it will continue to have two branches, one a
deformation of a
neighborhood of the point $p'_{i,j}$ in $L$ and the other a deformation of a
neighborhood of $p'_{i,j}$ in $X$. Similarly, a deformation $X_t$ of $X_0$
satisfying
condition iv) will have two branches near each point $s_{i,j}$,
deformations of the two
branches of $X_0$ at
$s_{i,j}$. In these terms, we make the further requirements that:

\

\item{v)}  In a neighborhood of each point $p'_{i,j}$, the branch of $X_t$
that is
a deformation of a neighborhood of $p'_{i,j}$ in $X$ has  contact of order
$i$ with
$L$ at $p'_{i,j}$
\smallskip
\item{vi)} In a neighborhood of each point $s_{i,j}$, the branch of $X_t$
that is
a deformation of a neighborhood of $s_{i,j}$ in $X$ has a point of
contact of order
$i$ with
$L$.

\

\ni \underbar{Remarks}.
1. We are being colloquial here in the definition of the relaxed
Severi variety, using terms like ``nearby" and ``in a neighborhood of" each
point
$p'_{i,j}$ or
$s_{i,j}$. This is to avoid introducing yet more notation. The definition
may be made
precise, for example, by specifying in
$\P^2$ disjoint analytic neighborhoods $U_{i,j}$, $V_{i,j}$ and
$W_1,\ldots,W_{\d'}$ of
the points
$p'_{i,j}$ and
$s_{i,j}$, and the nodes $u_1,\ldots,u_{\d'}$ of $X$; or by considering
deformations
$\{X_t\}$ of
$X_0$ having nodes at deformations $\{u_i(t)\}$ of the points $u_i$, etc.

2. We remark again that the relaxed local Severi variety $W_\Lambda$
depends on the
choice of subset $\{r_{i,j}\} \subset \{r'_{i,j}\}$; there are thus ${\b'
\choose \b}$ such
varieties
$W_\Lambda$ in a neighborhood of $[X_0]$.

\

Note that we are at this point making no requirements about the
deformations $X_t$ in
a neighborhood of a point $s''_{i,j}$, even though we have seen that a
family of curves
$X_t$ in $\V$ tending to $X_0$ will have $i-1$ nodes tending to each point
$s''_{i,j}$
(hence the name ``relaxed"). Thus, in particular, a general point $[X_t]
\in W_\Lambda$ will
correspond to a curve $X_t$ with only $\d'' = \d - (I\b'' - |\b''|)$
nodes---in other words,
$W_\Lambda$ will be an open subset of the variety $V^{d,\d''}(\a,\b)$. In
fact, our strategy is
exactly this: to consider first the conditions a curve $X_t$ must satisfy
away from the
points $s''_{i,j}$ in order to belong to a component of $\V$ containing
$[X_0]$ in its
closure; and then secondly the conditions around the points
$s''_{i,j}$.

The point
is, the conditions on the curve $X_t$ at the points
$p'_{i,j}$ and
$s'_{i,j}$ and the nodes $u_1,\ldots,u_{\d'}$ of $X$ are all essentially
linear conditions,
and well behaved. Thus, omitting any requirements on the behavior of the
curves
$X_t$ around
$s''_{i,j}$ will result in a  parameter space
$W_\Lambda$ that is smooth with identifiable tangent space at the point
$[X_0]$. Once
we have described this space, we will then consider the map
$\phi_\Lambda$ from
$W_\Lambda$ to the product $\Delta$ of the deformation spaces of the tacnodes of
$X_0$ at the points
$s''_{i,j}$. In a neighborhood of $[X_0]$, the Severi variety $\V$ will be
the union,
over all $\Lambda$, of the closures of the inverse images
$\phi_\Lambda^{-1}(\Delta_{m-1}
\setminus \Delta_m)$. Once we have shown that for each $\Lambda$ the image
$\phi_\Lambda(W_\Lambda)$ satisfies the hypotheses of Lemma 4.3, then,
Theorem 1.3 will
follow.

The first step in carrying
out this plan is thus the identification of the tangent space to
$W_\Lambda$ at $[X_0]$
(from which it will follow that $W_\Lambda$ is indeed smooth at $[X_0]$, once we
estimate its dimension). This tangent space, viewed as a subspace of the
tangent space $T_{[X_0]}\P^N = H^0(X_0,
\O(d))$, is the subspace $H^0(X_0,
\I(d))$ determined by an ideal sheaf $\I \subset \O_{X_0}$,
which we will describe in the following Lemma.

To do this, we have to introduce some
local ideals. Specifically, for each $i$ and $j$ with $1 \le j \le \a'_i$,
we let $\I'_{i,j}
\subset \O_{X_0}$ be the sheaf of regular functions in a neighborhood of
$p'_{i,j} \in
X_0$ whose restriction to  $L \subset X_0$ vanishes to order $i$ at
$p'_{i,j}$ and whose
restriction to  $X \subset X_0$ vanishes to order
$2i$ at $p'_{i,j}$. Similarly, for each $i$ and $j$ with $1 \le j \le
\b_i$, we let $\I_{i,j}
\subset \O_{X_0}$ be the sheaf of regular functions in a neighborhood of
$s_{i,j} \in
X_0$ whose restriction to $L \subset X_0$ vanishes to order $i$ at
$s_{i,j}$ and whose
restriction to  $X \subset X_0$ vanishes to order
$2i-1$ at $s_{i,j}$. We have then the

\

\proclaim Lemma 4.9.  The variety $W_\Lambda$ is smooth at $[X_0]$, and its
tangent space is given by the linear series
$$
T_{[X_0]}W_\Lambda = H^0(X_0,
\I(d)) \; \subset \; T_{[X_0]}\P^N = H^0(X_0,
\O(d))
$$
where the ideal sheaf $\I$ is the product
$$
\I \; = \; \prod_{i=1}^{\d'}\gm_{u_i} \cdot \prod_{1 \le j \le
\a''_i}\gm_{p''_{i,j}}^i \cdot
\prod_{1 \le j \le \a'_i} \I'_{i,j} \cdot \prod_{1 \le j \le \b_i} \I_{i,j}
$$

\

\ni {\it Proof}. We will first show that the tangent space to $W_\Lambda$ at
$[X_0]$ is contained in the series $H^0(X_0,
\I(d))$. We will then argue that the ideal $\I$ imposes independent
conditions on the
series $H^0(X_0,
\O(d))$, so that we can calculate the dimension $h^0(X_0,
\I(d))$; comparing this with the actual dimension of $W_\Lambda$ we deduce the
smoothness of $W_\Lambda$ at $[X_0]$ and the identification of the tangent space
with $H^0(X_0,
\I(d))$.

For both parts, it will be useful to  introduce a partial
normalization
$\tilde X_0$ of
$X_0$: specifically, we let $\mu : \tilde X_0 \to X_0$ be the normalization
of $X_0$ at the
points
$p'_{i,j}$ and
$s_{i,j}$, and the nodes $u_1,\ldots,u_{\d'}$ of $X$, but not at the points
$s''_{i,j}$. Note that the
normalization $\tilde X$ of $X$ is actually a closed subscheme of $\tilde
X_0$. We will
abuse notation slightly and denote by $q'_{i,j}$ and $r_{i,j}$ the points
of $\tilde X \subset \tilde
X_0$ lying over $p'_{i,j}$ and $s_{i,j}$. We will also denote by $q''_{i,j}
\in \tilde X_0$ the
(unique) point of $\tilde X_0$ lying over $p''_{i,j}$.

Note that since
$\I$ is contained in the conductor ${\cal J} \subset \O_{X_0}$ of the map
$\tilde X_0 \to
X_0$, the pullback
$\mu^*\I$ is locally free on $\tilde X_0$; specifically, it is the sheaf
$$
\mu^*\I \; = \; \mu^*{\cal J} \otimes \O_{\tilde X_0}\left(-\sum i\cdot
q'_{i,j} -\sum i\cdot
q''_{i,j} -
\sum (i-1)\cdot r_{i,j}\right) .
$$
Also, because $\I$ is contained in the conductor ${\cal J}$ of the map
$\tilde X_0 \to
X_0$, the space of sections $H^0(\tilde X_0, \mu^*\I(d)) = \mu^*H^0(X_0,
\I(d))$.

\

Now, to prove the inclusion $T_{[X_0]}W_\Lambda \subset H^0(X_0,
\I(d))$, we observe that by conditions i), iii) and iv) in the definition of
$W_\Lambda$ any deformation of $X_0$ in $W_\Lambda$ arises from a
deformation of the
composite map $\overline\mu : \tilde X_0 \to X_0 \to \P^2$. Thus, any
tangent vector to
$W_\Lambda$ at $[X_0]$ must lie in the subspace
$$
H^0(X_0, {\cal J}(d)) \; \subset \; H^0(X_0, \O(d))
$$
which we may identify in turn with $H^0(\tilde X_0, \mu^*{\cal J}(d))$.
Now, as we observed in
section 2 above, we may further identify $\mu^*{\cal J}(d)$ with the normal
sheaf
$N = N_{\overline\mu}$ of the map $\overline\mu$;  and in terms of these
identifications
conditions ii), v) and vi) of the definition of $W_\Lambda$ amount to the
assertion that the
tangent space to $W_\Lambda$ at $[X_0]$ satisfies
$$
\eqalign{T_{[X_0]}W_\Lambda \; &\subset \; H^0(\tilde X_0, N(-\sum i\cdot
q'_{i,j} -\sum i\cdot
q''_{i,j} -
\sum (i-1)\cdot r_{i,j})) \cr
&= H^0(\tilde X_0, \mu^*{\cal J}(d))(-\sum i\cdot q'_{i,j} -\sum i\cdot
q''_{i,j} -
\sum (i-1)\cdot r_{i,j})) \cr
&= H^0(\tilde X_0, \mu^*\I(d)) \cr
&= H^0(X_0, \I(d)) \cr}
$$

\

For the second part, to estimate of the dimension $h^0(X_0,
\I(d))$, we will equate this with the dimension $h^0(\tilde X_0,
\mu^*\I(d))$ of the space of
sections of the pullback, and apply Riemann-Roch on $\tilde X_0$. A key
fact is that the line
bundle $\mu^*\I(d)$ is nonspecial. In the sequel we will need as well the
fact that for $p \in L$
general the bundle $\mu^*\I(d)(-p)$ is nonspecial as well; we will state
these as the

\

\proclaim Lemma 4.10. For $p \in L$ general,
$$
h^1(\tilde X_0, \mu^*\I(d)) \; = \; h^1(\tilde X_0, \mu^*\I(d)(-p)) \; = \;
0 \, .
$$

\

\ni {\it Proof}.
First, note that the dualizing sheaf of $\tilde X_0$ is given by
$$
\omega_{\tilde X_0} \; = \; \mu^*({\cal J}(d-3)) .
$$
Thus, we may write
$$
\mu^*\I(d) \; = \; \mu^*\O(3) \otimes \omega_{\tilde X_0}\left(-\sum i\cdot
q'_{i,j} -\sum
i\cdot q''_{i,j} -
\sum (i-1)\cdot r_{i,j}\right) .
$$
and correspondingly
$$
\omega_{\tilde X_0} \otimes (\mu^*\I(d))^{-1} \; = \; \mu^*\O(-3)\left(\sum
i\cdot q'_{i,j} +\sum
i\cdot q''_{i,j} +
\sum (i-1)\cdot r_{i,j}\right) \, .
$$
Now, the restriction of this line bundle to $\tilde X$  has degree
$$
\eqalign{\deg\left(\omega_{\tilde X_0} \otimes (\mu^*\I(d))^{-1} \otimes
\O_{\tilde X}\right) \;
&=
\; -3(d-1) + I\a' + I\b \cr
&\le
\; -3(d-1) + I\a' + I\b' \cr
&= \; -2(d-1) \cr
 &< \; 0 \cr}
$$
and so every global  section $\s \in H^0(\omega_{\tilde X_0} \otimes
(\mu^*\I(d))^{-1})$ must
vanish identically on $\tilde X$. The restriction of $\s$ to the component
$\tilde L$ of $\tilde
X_0$ lying over $L$ is then a section of the bundle
$$
\omega_{\tilde X_0} \otimes (\mu^*\I(d))^{-1} \otimes \I_{\tilde X} \otimes
\O_{\tilde L}
$$
and since the restriction $\I_{\tilde X} \otimes \O_{\tilde L}$ to $\tilde
L$ of the ideal sheaf of
$\tilde X$ has degree $-I\b''$, we have
$$
\eqalign{\deg\left(\omega_{\tilde X_0} \otimes (\mu^*\I(d))^{-1} \otimes \I_{\tilde X} \otimes
\O_{\tilde L}\right) \; &= \; I\a'' - I\b'' -3 \cr
&= \;  I\a - I\a' - (I\b' - I\b) -3  \cr
&= \;  (I\a+I\b) - (I\a' + I\b') -3 \cr
&= \; d - (d-1) - 3 \cr
&= \; -2 \, .\cr}
$$
Thus $\s$ must vanish identically on $\tilde L$ as well, and hence
$$
h^1(\tilde X_0, \mu^*\I(d)) \; = \; h^0(\tilde X_0, \omega_{\tilde X_0}
\otimes (\mu^*\I(d))^{-1})
\; = \; 0 \, .
$$
Finally, if we had started with $\mu^*\I(d)(-p)$ in place of $\mu^*\I(d)$,
we would have
wound up with
$$
\deg\left(\omega_{\tilde X_0} \otimes (\mu^*\I(d)(-p))^{-1} \otimes
\I_{\tilde X} \otimes
\O_{\tilde L}\right) \; = \; -1
$$
and we would conclude as before that the line bundle $\mu^*\I(d)(-p)$ is
nonspecial.
\qed

\

We may now apply the Riemann-Roch formula on $\tilde X_0$ to complete the
proof of Lemma
4.9. By Lemma 4.10, we have
$$
\eqalign{h^0(X_0, \I(d)) \; &= \; h^0(\tilde X_0, \mu^*\I(d)) \cr
&= \deg(\mu^*\I(d)) - p_a(\tilde X_0) + 1   \cr}
$$
Since
$$
\eqalign{\deg(\mu^*\I(d)) \; &= \; \deg\left(\mu^*\O(3) \otimes
\omega_{\tilde X_0}(-\sum
i\cdot q'_{i,j} -\sum i\cdot q''_{i,j} -
\sum (i-1)\cdot r_{i,j})\right) \cr
&= \; 3d +2p_a(\tilde X_0) - 2 -I\a -(I\b - |\b|) \cr}
$$
we can rewrite this as
$$
h^0(X_0, \I(d)) \; = 3d + p_a(\tilde X_0) - 1 -I\a -(I\b - |\b|) \, .
$$
Now, the arithmetic genus of $\tilde X_0$ is simply the genus of a general
member $X_t$ of
$W$; thus
$$
p_a(\tilde X_0) \; = \; {d-1 \choose 2} -\d'' .
$$
Equivalently, we could arrive at this by observing that the arithmetic
genus of $\tilde X_0$ is
simply the genus of $\tilde X$, plus the degree of the intersection of
$\tilde X
\subset
\tilde X_0$ with the component
$\tilde L$ of
$\tilde X_0$ lying over the line $L \subset \P^2$, minus 1; thus
$$
\eqalign{p_a(\tilde X_0) \; & = \; {d-2 \choose 2} -\d' + I\b'' - 1 \cr
&= \; {d-2 \choose 2} -(\d + d - 1 -|\b'-\b|)  + I\b'' - 1 \cr
&=  \; {d-1 \choose 2} -\d + I\b'' - |\b''|  \cr
&= \; {d-1 \choose 2} -\d'' \cr} .
$$
Either way, we have
$$
\eqalign{\dim(T_{[X_0]}W_\Lambda) \; &\le \; h^0(X_0, \I(d)) \cr
 &= 3d + p_a(\tilde X_0) - 1
-I\a -(I\b - |\b|)
\cr &= 3d + {d-1 \choose 2} -\d''  - 1 -I\a -(I\b - |\b|)  \cr
&= {d+2 \choose 2} - 1 - \d''   -I\a -(I\b - |\b|) \cr
&= \dim(V^{d,\d''}(\a,\b)) \cr
&= \dim(W_\Lambda) \, ; \cr}
$$
so $W_\Lambda$ must be smooth, with tangent space equal to $H^0(X_0, \I(d))$.
\qed

\

All that remains to complete the proof of Proposition 4.8 is to consider the map
$\phi_\Lambda$ from
$W_\Lambda$ to the product $\Delta$ of the deformation spaces of the tacnodes of
$X_0$ at the points
$s''_{i,j}$. Note first that the image of $\phi_\Lambda$ contains the
locus $\Delta_m \subset \Delta$, and that the inverse image $W_0 =
(\phi_\Lambda)^{-1}(\Delta_m)$ of this locus  is simply the set of points
$[X_t] \in W_\Lambda$ corresponding to curves $X_t$ containing $L$. This
has codimension at
most one in $W_\Lambda$: if we choose any point $p \in L$ not among the
points $p_{ij}$ or
$s_{i,j}$, then any curve $X_t$ in $W_\Lambda$ containing $p$ will have a
total of
$$
I\a + I\b + 1 \; = \; d+1
$$
points of intersection with $L$, and so will contain $L$.

\

On the other hand, under the differential
$$
d\phi_\Lambda \; : \; T_{[X_0]}W_\Lambda \; \la \; T_0\Delta
$$
of the map $\phi_\Lambda$, the inverse image of the tangent space
$T_0\Delta_m \subset
T_0\Delta$ is simply the subspace  of $H^0(X_0,\I(d))$ of sections
vanishing on $L$. Now, the
restriction of $\I(d)$ to $L$ has degree $d - I\a - I\b = 0$---that is, it
is trivial---and by
Lemma 4.10 not every global section of $\I(d)$ vanishes on $L$. Thus
$$
\dim((d\phi_\Lambda)^{-1}(T_0\Delta_m)) \; = \; \dim(\phi^{-1}(\Delta_m)
\; = \;
\dim(W_\Lambda) -1 .
$$
We may conclude that the image of $\phi_\Lambda$ is smooth of dimension
$$
\dim(\phi_\Lambda(W_\Lambda)) \; = \; \dim(\Delta_m) + 1
$$
with tangent space the image of $d\phi_\Lambda$. Since again the linear system
$H^0(X_0,\I(d))$ has no base points on $L$, the image $W =
\phi_\Lambda(W_\Lambda)$
satisfies the hypotheses of Lemma 4.3. Thus  we may apply Lemma 4.3 to conclude that the
closure of the inverse image $\phi_\Lambda^{-1}(\Delta_{m-1}
\setminus \Delta_m)$ will have $\k = I^{\b'-\b}/\lcm(\b' - \b)$ reduced
branches, each having
intersection multiplicity $\lcm(\b'-\b)$ with $W_0$ and hence with the
hyperplane $H_p$.
Since in a neighborhood of
$[X_0]$ the Severi variety
$$
\V \; = \; \bigcup_\Lambda \overline{\phi_\Lambda^{-1}(\Delta_{m-1}
\setminus \Delta_m)}
$$
we conclude finally that near $[X_0]$, $\V$ will have ${\b' \choose
\b}I^{\b'-\b}/\lcm(\b' -
\b)$  branches, each of which will  have intersection multiplicity
$\lcm(\b'-\b)$ with $H_p$
along
$V'$.  This completes the proof of Proposition 4.8 and thereby of Theorem 1.3
\qed

\vfill\eject

\

\

\cl{\sc References}

\

\ni [CH] L. Caporaso and  J. Harris, {\it Parameter spaces for curves on
surfaces and enumeration of rational curves.}
\ (1995)  preprint

\ni [DH] \ S. Diaz and J. Harris, {\it Geometry of the Severi variety}. \
Trans.A.M.S. 309 (1988),
pp. 1-34

\ni [H] \  J. Harris, {\it On the Severi problem}.  \ Invent. Math. 84
(1986), pp. 445-461.

\ni [K] \  J.Kollar, {\it Rational curves on algebraic varieties}.
Springer (1996) .

\ni [KM] M. Kontsevich and Y. Manin,  \ {\it Gromov-Witten classes, quantum
cohomology
and enumerative geometry.}  preprint.

\ni [P] R. Pandharipande, {\it Intersections of $\Q$-divisors on
Kontsevich's moduli space}

\ni {\it $\overline{M}_{0,n}(\P^r, d)$ and enumerative geometry}. \ (1995),
preprint

\ni [R] Z. Ran, {\it Enumerative geometry of singular plane curves}.
Invent. Math. 97 (1989), pp. 447-465.

\ni [RT] Y. Ruan and G. Tian, {\it A
mathematical theory of quantum cohomology}, J. Diff. Geom. 42 No. 2 (1995)
pp. 259-367).

\ni [Va] I. Vainsencher,  {\it Counting divisors with prescribed
singularities.} Trans AMS 267 (1981), 399-422.

\ni [V] R. Vakil, {\it Curves on  rational ruled surfaces.}
(1996) preprint.

\end